\documentclass[10pt]{article}
\usepackage[latin1]{inputenc}
\usepackage{amsmath,amssymb}
\usepackage{color}
\usepackage{multicol}
\usepackage{latexsym,epsfig}
\usepackage{amssymb,amsmath,amsfonts,amscd}
\usepackage{exscale}
\usepackage{ifthen}
\usepackage{longtable}
\usepackage{subfigure}
\usepackage{cite}
\usepackage{lscape}

\title{Counting for some convergent groups}

\newtheorem{theo}{\textbf{Theorem}\ }
[section]
\newtheorem{lemma}[theo]{\textbf{Lemma}\ }
\newtheorem{coro}[theo]{Corollary\ }

\newtheorem{prop}[theo]{\textbf{Proposition}\ }

\newtheorem{property}[theo]{\textbf{Property}\ }

\newtheorem{notations}[theo]{Notations\ }

\newtheorem{remark}[theo]{\textbf{Remark}\ }

\newtheorem{fact}[theo]{\textbf{Fact}\ }

\begin{document}     

\maketitle

 \centerline{  Marc Peign\'e $^($\footnote{Marc
Peign\'e, LMPT, UMR 7350, Facult\'e des Sciences et Techniques,
Parc de Grandmont, 37200 Tours, France -- email :
peigne@univ-tours.fr}$^)$, Samuel Tapie $^($\footnote{Samuel Tapie,
Laboratoire Jean Leray 
2 rue de la Houssini\`ere - BP92208 
44322 NANTES CEDEX 3
France
 -- email :
samuel.tapie@univ-nantes.fr}$^)$ $\&$ Pierre Vidotto
$^($\footnote{Pierre Vidotto,
Laboratoire Jean Leray 
2 rue de la Houssini\`ere - BP92208 
44322 NANTES CEDEX 3
France
 -- email :
pierre.vidoto@univ-nantes.fr}$^)$}

 \vspace{0.5cm}

\noindent {\bf Abstract.}
We present examples of geometrically finite manifolds  with pinched negative curvature, whose geodesic flow has infinite non-ergodic Bowen-Margulis measure and whose Poincar\'e series converges at the critical exponent $\delta_\Gamma$. We obtain an explicit asymptotic for their orbital growth function.  Namely, for any $\alpha \in ]1, 2[ $ and any slowly varying function $L : \mathbb R\to (0, +\infty)$, we construct $N$-dimensional Hadamard manifolds $(X, g)$ of negative and pinched curvature, whose group of oriented isometries admits  convergent geometrically finite subgroups $\Gamma$  such that, as $R\to +\infty$,
$$
N_\Gamma(R):= \#\left\{\gamma\in \Gamma \; ; \; d(o, \gamma \cdot o)\leq R\right\} \sim C_\Gamma \frac{L(R)}{R^\alpha} \ e^{\delta_\Gamma R},
$$
for some constant $C_\Gamma >0$.

\vspace{3mm}


\vspace{3mm}

\section{Introduction}

We fix  $N\geq 2$ and consider  a $N$-dimensional Hadamard manifold $X$ of negative, pinched curvature $-B^2 \leq  K_X \leq -A^2<0$. Without loss of generality, we may assume $A\leq 1 \leq B$. Let $\Gamma$ be a {\em Kleinian group} of $X$,  i.e. a discrete, torsionless group of isometries of  $X$, with quotient $\bar X= \Gamma \backslash X$.

This paper is concerned with the  fine asymptotic properties of the {\em orbital function} : 
$$v_{\Gamma}({\bf x},{\bf y};R):=\sharp\{\gamma\in\Gamma\slash d({\bf x},\gamma\cdot{\bf y})\leq R\}$$
for ${\bf x},{\bf y} \in X$, which has been the subject of many investigations since Margulis' \cite{Ma} (see also Roblin's book \cite{Ro}). First, a simple   invariant   is its {\em exponential growth rate}
$$\delta_{\Gamma}=\limsup_{R\to\infty}\frac{1}{R}\log(v_{\Gamma}({\bf x},{\bf y};R)).$$
The exponent $\delta_\Gamma$ coincides also with the {\em exponent of convergence} of the {\it Poincar\'e series} associated with $\Gamma$ :
$$P_{\Gamma}({\bf x}, {\bf y}, s):=\sum_{\gamma\in\Gamma}e^{-sd({\bf x},\gamma\cdot{\bf y})}, \qquad {\bf x}, {\bf y}\in X. $$
Thus, it is called the {\it Poincar\'e exponent} of $\Gamma$ $\delta_\Gamma$. It coincides with the topological entropy of the geodesic flow $(\phi_t)_{t\in \mathbb R}$ on the unit tangent bundle of $\bar X$, restricted to its non-wandering set. It equals also  the Hausdorff dimension of the {\it radial limit set} $\Lambda(\Gamma)^{rad}$ of $\Gamma$ with respect to some natural metric on the  boundary at infinity $\partial X$ of $X$. Recall that  any orbit $\Gamma\cdot{\bf x}$ accumulates on  a  closed subset  $\Lambda(\Gamma)$ of the  geometric boundary $\partial X$ of $X$, called the {\em limit set} of $\Gamma$; this set contains 1, 2 or infinitely many points and one says that $\Gamma$ is non elementary when $\Lambda_\Gamma$ is infinite. 
 A point $x \in\Lambda_{\Gamma}$ is said to be   {\em radial} when it is approached by orbit points in some $M$-neighborhood of any given ray   issued from $x$, for some $M>0$). 
 
The group $\Gamma$ is said to be {\em convergent} if $P_{\Gamma}({\bf x}, {\bf y}, \delta_{\Gamma})<\infty$, and {\em divergent } otherwise. Divergence can  also be understood in terms of  dynamics as, by Hopf-Tsuju-Sullivan theorem, it is equivalent to ergodicity and total conservativity of the geodesic flow with respect to the Bowen-Margulis measure $m_{\Gamma}$  on the  non wandering set of $(\phi_t)_{t\in \mathbb R}$  in the unit tangent bundle $T^1\bar X$ (see again \cite{Ro} for a complete account and a definition of  $m_{\Gamma}$ and for a proof of this equivalence).

The more general statement concerning the asymptotic behavior of 
$v_{\Gamma}({\bf x},{\bf y};R)$ is due to Th. Roblin: if $\Gamma$ is a     non elementary, discrete subgroup of isometries of $X$ with {\em non-arithmetic length spectrum}\footnote{It means that the set 
$\mathcal{L} (\bar X)=\{  \ell (\gamma)\; ;\; \gamma\in\Gamma\}$  of lengths of  closed geodesics of $\bar X = \Gamma \backslash X$ is not contained in a discrete subgroup of $\mathbb R$}, then    $\delta_{\Gamma}$ is a true limit and it holds, as $R\to +\infty$,
\begin{enumerate}
\item [(i)]  if $\Vert m_{\Gamma}||=\infty$ then $v_{\Gamma}({\bf x},{\bf y};R)=o(e^{\delta_{\Gamma}R})$,
\item [(ii)]  if $\Vert  m_{\Gamma}\Vert <\infty$, then 
$v_{\Gamma}({\bf x},{\bf y};R)\sim{||\mu_{\bf x}||. ||\mu_{\bf y}||\over \delta_{\Gamma}\Vert m_{\Gamma}||}e^{\delta_{\Gamma}R},$
\end{enumerate}
 where $(\mu_{\bf x})_{{\bf x} \in X}$ denotes the family of Patterson $\delta_\Gamma$-conformal densities of $\Gamma$, and  $m_{\Gamma}$ the Bowen-Margulis measure on  $T^1\bar X$.
Let us emphasize that in the second  case, the group $\Gamma$ is always divergent while in the first  one it can be   convergent.

In this paper, we  aim to investigate, for a particular class of groups $\Gamma$,  the asymptotic behavior of the function $v_{\Gamma}({\bf x},{\bf y};R)$ when $\Gamma$ is convergent. 
As far as we know, the only precise asymptotic for the orbital  function of   convergent groups   holds  for groups $\Gamma$  which are  normal subgroups  $\Gamma \lhd \Gamma_0$  of a co-compact group $\Gamma_0$ for which the quotient $\Gamma_0\slash \Gamma$ is isometric up to a finite factor to the lattice  $\mathbb Z^k$ for some $k\geq 3$  \cite{PS}.  The corresponding quotient manifold has infinite Bowen-Margulis measure; in fact, $m_{\Gamma}$ is invariant under the action of the group of isometries of  $\bar X$ which contains subgroups $\simeq \mathbb Z^k$.

 The finiteness of $m_{\Gamma}$ is not easy  to establish excepted in  the case of {\it geometrically finite groups} where there exists a precise criteria. Recall that
 $\Gamma $  (or the quotient manifold $\bar X$)   is said to be geometrically finite if  its limit set $\Lambda(\Gamma)$  decomposes in the  radial limit set and the $\Gamma$-orbit of finitely many  {\em bounded parabolic points} $x_1, \ldots, x_\ell$,  fixed  respectively by some parabolic subgroups $P_i, 1\leq i \leq \ell$,  acting co-compactly on
  $\partial X  \setminus \{x_i\}$; for a complete  description of geometrical finiteness in variable negative curvature   see   \cite{Bow}. Finite-volume  manifolds $\bar X$ (possibly non compact)  are particular cases of geometrically finite manifolds; in contrast, the manifolds considered in \cite{PS} are not     geometrically finite.

For geometrically finite groups, the orbital functions 
$v_{P_i}$ of the parabolic subgroups $P_i, 1\leq i\leq \ell$,
contain  the relevant information about the metric   inside  the cusps, which in turn may imply   $ m_{\Gamma}  $ to be finite or infinite.  
On the one hand, it  is proved  in  \cite{DOP} that 
the divergence of the parabolic subgroups $P\subset \Gamma$ implies $\delta_{P } < \delta_\Gamma$, which in turn yields that 
 $\Gamma$ is divergent and  $\Vert m_{\Gamma}||<\infty$. 
 On the other hand there exist geometrically finite groups with  parabolic subgroups $P $ satisfying  $\delta_{P} = \delta_\Gamma$: 
  we call such groups {\em exotic}  and   say that the parabolic subgroup  $P $ (or the corresponding cusps $\mathcal C $)  is {\em dominant} when $\delta_{P} = \delta_\Gamma$.   Let us emphasize that  dominant parabolic subgroups of exotic geometrically finite groups  $\Gamma$ are necessarily convergent. However, the group $\Gamma$ itself may as well be convergent or divergent; we refer to \cite{DOP} and  \cite{P} for explicit constructions of such groups.

In this paper, we   consider a  Schottky product  $\Gamma$  of  elementary subgroups 
$ \Gamma_1, \ldots,   \Gamma_{p+q},$
   of isometries of   $X$ (see  $\S 3$ for the definition) with $ p+q \geq 3$. 
Such a group  is geometrically finite. We   assume that $\Gamma$ is convergent;  thus, by \cite{DOP},  it is exotic  and  possesses  factors (say $\Gamma_1,\ldots, \Gamma_p, p \geq 1$) which are  dominant  parabolic subgroups  of $\Gamma$. We assume that, up to the dominant factor $e^{\delta_\Gamma R}$, the orbital functions $v_{\Gamma_j}({\bf x}, {\bf y}, \cdot)$ of these groups  satisfy some asymptotic condition of polynomial decay at infinity.  More precisely we have the

\begin{theo}\label{theo:Comptage}
Fix  $p, q \in \mathbb  N$ such that $p \geq 1, p+q \geq 2$ and let $\Gamma$ be a Schottky product of  elementary subgroups 
$\Gamma_{1}, \Gamma_2\ldots, \Gamma_{p+q}$ of isometries of  a pinched negatively curved manifold $X$.
Assume that the metric $g $ on  $X$ satisfies the following assumptions.
\vspace{1mm}

 $  {\bf H_1.}$ The group $\Gamma$ is convergent  with Poincar\'e exponent  $\delta_\Gamma=\delta$.

$  {\bf H_2.}$    There exist $\alpha\in ]1, 2[$, a slowly varying  function  $L$\
$^($\footnote{A function $L(t)$ is said to be ''slowly varying''   it is positive, measurable and $L(\lambda t)/L(t)\to  1$  as $t\to +\infty$ for every  $\lambda>0$.}$^)$ and  strictly positive constants   $c_1, \ldots, c_p$    such that, for any $1\leq j \leq p$ and  $\Delta >0,$
\begin{equation}\label{influent}
 \lim_{R\to +\infty} {R^{\alpha} \over L(R)}
 \sum _{ \stackrel{\gamma \in \Gamma_j}{R\leq d({\bf o}, \gamma\cdot {\bf o})<R+\Delta}} 
 e^{-\delta d({\bf o}, \gamma\cdot {\bf o})}  =c_j \Delta.
 \end{equation}

$ {\bf H_3.}$ {\it For any $p+1\leq j\leq p+q$ and   $\Delta >0,$
$$
 \lim_{R\to +\infty} {R^{\alpha} \over L(R)}
 \sum _{ \stackrel{\gamma \in \Gamma_j}{R\leq d({\bf o}, \gamma \cdot {\bf o})<R+\Delta}} 
 e^{-\delta d({\bf o}, \gamma\cdot {\bf o})}  =0.  
$$}

\noindent Then,   there exists  a constant $C_\Gamma>0$ such that, as $R\to +\infty$, 
$$
\sharp\{\gamma\in \Gamma\mid d({\bf o}, \gamma \cdot {\bf o})\leq R \}
\quad 
\sim
\quad 
C_\Gamma\  {L(R)\over  R^{\alpha}}\  e^{\delta R}.
$$
\end{theo}

The importance of the convergence hypothesis {\bf H$_1$} in the previous theorem is illustrated by the following result, previous work of one of the authors. 

\begin{theo}[\cite{V}, Theorem C]
Let $\Gamma$ be a Schottky product of  $p+q\geq 2$ elementary subgroups 
$\Gamma_{1}, \Gamma_2\ldots, \Gamma_{p+q}$ of isometries of  a pinched negatively curved manifold $X$. Assume that  $p \geq 1$ and 

$\bullet$  $\Gamma$ is \emph{divergent}  and  $\delta_\Gamma = \delta$,
 
$\bullet$ Hypotheses {\bf H$_2$}, {\bf H$_3$} hold. 

\noindent Then, there exists $C_\Gamma>0$ such that, as $R\to +\infty$,
$$
\sharp\{\gamma\in \Gamma\mid d({\bf o}, \gamma \cdot {\bf o})\leq R \}
\quad 
\sim
\quad 
C_\Gamma\  {e^{\delta R}\over  R^{2-\alpha}L(R)}.
$$
\end{theo}

The difference with the equivalent of Theorem \ref{theo:Comptage} may surprise, since it is possible to vary smoothly the Riemannian metric $g_{\alpha, L}$ from a divergent to a convergent case, preserving hypotheses {\bf H$_2$} and {\bf H$_3$}, cf \cite{P}. Nevertheless, the proof of our Theorem \ref{theo:Comptage} will illustrate the reasons of this difference. 
For groups $\Gamma = \Gamma_1*...*\Gamma_{p+q}$ satisfying {\bf H$_2$} and {\bf H$_3$}, the counting estimate only depends on elements of the form 
$\gamma =a_1 \cdots a_k$, whith $a_i\in \Gamma_1\cup\ldots \cup \Gamma_p$ and where $a_i$ and  $a_{i+1}, 1\leq i <k$,  do not belong to  the same $\Gamma_j$. In the divergent case (see the proof of Theorem C in\cite{V}), the asymptotic of $\{\gamma\in \Gamma\mid d({\bf o}, \gamma \cdot {\bf o})\leq R \}$ as $R\to +\infty$ only depends on the $\gamma = a_1\cdots a_k$ with $k>> R$. On the opposite,  in the convergent case,   the dominant parabolic factors $\Gamma_1, \ldots, \Gamma_p$ are ``heavy'' and the  asymptotic  of the orbital function of $\Gamma$ comes from the $\gamma = a_1\cdots a_k$ with $k$ bounded independently of $R$; the number of such isometries $\gamma$ with $d({\bf o}, \gamma \cdot {\bf o}) \leq R$ is comparable to $\displaystyle {L(R)\over  R^{\alpha}}\  e^{\delta R}$. By a straightforward adaptation of  Proposition \ref{mkequi}, this last estimate  remains valid in the divergent case; nevertheless,  the fact that $\Gamma$ is divergent implies that the contribution of these isometries is negligible.

\begin{remark}
The condition $\alpha >1$ assures that the parabolic groups $\Gamma_1, \ldots, \Gamma_p$ are convergent. The additive condition $\alpha<2$ is used in Proposition \ref{mkmaj} to obtain a uniform upper bound for the power $\widetilde P^k, k\geq 1$ of some operator $\widetilde P$  introduced in Section  \ref{sec:Counting}; the proof of this Proposition relies on a previous work of one of the authors \cite{V} and is not valid for  greater values of $\alpha$. The  analogous of our Theorem \ref{theo:Comptage} when $\alpha \geq 2$ remains open.
\end{remark}

The article is organized as follows. In the next section, we recall some backgrounds on negatively curved manifolds, and we construct examples of metrics for which the hypotheses of Theorem \ref{theo:Comptage} are satisfied. In section \ref{sec:Schottky}, we present Schottky groups and the coding which we  use to express our geometric problem in  terms of  sub-shift of finite type on a  countable alphabet. In section \ref{sec:Ruelle}, we introduce the Ruelle operator for this sub-shift; this is the key analytical tool which is used. Eventually, section \ref{sec:Counting} is devoted to the proof of Theorem \ref{theo:Comptage}.


\section{Geometry of negatively curved manifolds}\label{sec:GeoNeg}

\subsection{Generalities}

 In the sequel, we fix $N \geq 2$ and consider a $N$-dimensional complete connected Riemannian manifold $X$ with metric $g$ whose sectional curvatures satisfy : $-B^2\leq K_X\leq -A^2<0$ for fixed constants $A$ and $B$; the metric $g$ we consider in this paper    be obtained by perturbation of a  hyperbolic one and the curvature    equal  $-1$ on large subsets of $ X$, thus we  assume 
$0<A\leq 1\leq B$.  We denote $d$ the distance on $X$ induced by the metric $g$.

Let $\partial X$ be the boundary at infinity of $X$ and let us fix an origin ${\bf o} \in X$.
The family of functions  $ \left({\bf y} \mapsto 
 d({\bf o},{\bf x})-d({\bf x}, \bf y )\right)_{ {\bf x} \in X}$   converges uniformly on compact sets to the  {\em Busemann function}  $ \mathcal B_{x}( {\bf o}, \cdot ) $  
for  ${\bf x}\to x\in \partial X$. 
The {\em horoballs} $\mathcal{H}_{x}$    and the  {\em horospheres} $\partial \mathcal{H}_{x}$   centered at $x$  are   respectively the sup-level sets and  the  level sets  of the function $ \mathcal B_{x}( {\bf o}, \cdot ) $.
For any $t\in\mathbb R$,   we set 
  $\mathcal{H}_{x}(t):=\{{\bf y}\slash  \mathcal B_{x}( {\bf o}, {\bf y} )   \geq t \}$
and  $\partial\mathcal{H}_{x}(t):=\{{\bf y}\slash \mathcal B_{x}( {\bf o}, {\bf y} ) =  t\}$; the parameter  $t = \mathcal B_{x}( {\bf o}, {\bf y}) -  \mathcal B_{x}( {\bf o},{\bf x})$ is the {\em height} of ${\bf y}$  with respect to $x$. When no confusion is possible, we   omit  the index $x \in \partial X$ denoting the center of the horoball. 
Recall that the Busemann function satisfies the fundamental cocycle relation: for any $x \in \partial X$ and  any $\bf x, y, z$ in $X$ 
$$\mathcal B_{x}(  {\bf x}, {\bf z}) = B_{x}(  {\bf x}, {\bf y}) + B_{x}(  {\bf y}, {\bf z}).
$$
The Gromov product between    $x, y \in \partial X \cong \partial X$,   $x \neq y$, is  defined as 
$$(x\vert y)_{{\bf o}} = {\mathcal B_x({\bf o}, {\bf z})+\mathcal B_y({\bf o}, {\bf z})\over 2}$$
 where ${\bf z}$ is any point on the geodesic $(x, y)$ joining $x$ to $y$. 
By \cite{Bou}), the expression
$$ D(x, y)= e^{ -A(x\vert y)_{{\bf o}} }$$
defines  a distance on $\partial X $   satisfying the following property: for any $\gamma \in \Gamma$  
$$
D(\gamma  \cdot x, \gamma \cdot y)= e^{-{A\over 2}  \mathcal B_x(\gamma^{-1}\cdot {\bf o}, {\bf o})} e^{-{A\over 2}  \mathcal B_y(\gamma^{-1}\cdot {\bf o}, {\bf o})} D(x, y).
$$
In other words, the isometry $\gamma$ acts on $(\partial X, D)$ as a conformal transformation with coefficient of conformality 
$ 
\vert \gamma'(x)\vert_{\bf o} = e^{-A  \mathcal B_x(\gamma^{-1}\cdot {\bf o}, {\bf o})}
$ 
at  $x$ and satisfies the following equality
\begin{equation} \label{TAF1}
D(\gamma  \cdot x, \gamma \cdot y)= \sqrt{\vert \gamma'(x)\vert_{\bf o} \vert  \gamma'(y)\vert_{\bf o} }  D(x, y).
\end{equation}  
The function   $x\mapsto \mathcal B_x(\gamma^{-1}\cdot {\bf o}, {\bf o})$ plays a central role to describe the action of the isometry $\gamma$ on  the boundary at infinity $\partial X$. From now on, we   denote it $b(\gamma, \cdot )$  and notice that it satisfies the following ``cocycle property'': for any  isometries $\gamma_1, \gamma_2$  of $X$ and any $x \in  \partial X$
\begin{equation} \label{cocycleb}
b(\gamma_1 \gamma_2, x) = 
b(\gamma_1,  \gamma_2 \cdot x)+ b( \gamma_2, x).
\end{equation} 
In order to describe the action  on $\partial X$ of the isometries of $(X, g)$, it is useful to control precisely the behavior of the sequences $\vert (\gamma^n)'(x)\vert_{\bf o}$; the following fact provides a useful estimation of these quantities.
\begin{fact}\label{lienentrebetd}
(1)
For any hyperbolic isometry $h$ with repulsive and attractive  fixed point $\displaystyle x_h^-=\lim_{n\to +\infty} h^{-n}\cdot{\bf o}$ and  $ \displaystyle x_h^+=\lim_{n\to +\infty} h^{n}\cdot{\bf o}$ respectively, it holds
$$
 b(h^{\pm n}, x) =d(o, h^{\pm n}\cdot o)-2(x_h^{\pm}\vert  x)_o+\epsilon_x(n) 
$$
with $\displaystyle \lim_{n\to+\infty}\epsilon_x(n)=0$, the convergence being uniform on the compact sets of $\partial X\setminus\{x_h^{\mp}\}$.

(2)
For any parabolic group  $\mathcal P$ with   fixed point $\displaystyle x_{\mathcal P}:=\lim_{\stackrel{p \in {\mathcal P}}{p\to +\infty}} p\cdot o$, it holds
$$
  b(p, x)=d(o,  p\cdot o)-2(x_\mathcal P\vert  x)_o+\epsilon_x(p) 
$$
with $\displaystyle \lim_{\stackrel{p \in \mathcal P}{p\to +\infty}} \epsilon_x(p)=0$, the convergence being uniform on the compact sets of $\partial X\setminus\{x_{\mathcal P} \}$.
\end{fact}

\subsection{On the existence of convergent parabolic groups}

 In this section, we recall briefly the construction   presented in \cite{P} of convergent parabolic groups satisfying  condition (\ref{influent}), up to a bounded term; we refer to \cite{P} for the details.

We consider on $\mathbb R^{N-1}\times \mathbb R$   a
Riemannian metric   of the form
$g=T^2(t){\rm d}x_{ }^2+{\rm d}t ^2$ at point ${\bf x} = (x, t)$
where ${\rm d}x_{ }^2$ is a fixed   euclidean  metric 
on  $\mathbb R^{N-1}$ and  $T: \mathbb R\to \mathbb R^{*+}$ is
a $C^{\infty}$   non-increasing function.
The group of isometries of  $g$ contains the isometries of 
$\mathbb R^{N-1}\times \mathbb R$   fixing the last  coordinate. The sectional 
curvature 
at  $\displaystyle {\bf x}= (x,t) $ equals $\displaystyle K_g(t)=-\frac 
{\ \ T''(t)}{T(t)}$ on any   plane  
$\displaystyle \Big\langle \frac {\partial }{\partial X_i}, \frac 
{\partial}{\partial t}\Big\rangle, 1\leq i\leq N-1$, 
and  $-K_g^2(t)$ on any plane  
$\displaystyle \Big\langle \frac {\partial }{\partial X_i}, \frac 
{\partial}{\partial X_j}\Big\rangle, 1\leq i<j\leq N-1$. Note that  $g$ has negative curvature if and only if $T$ is 
convex; when $T(t)= e^{-t}$, one obtains a 
model of the hyperbolic space   of constant curvature $-1$.

   It is convenient to consider the 
  non-decreasing function
\begin{eqnarray}\label{eq:FonctionU}
u :\left\{\begin{array}{lll}
 \mathbb R^{*+}&\to &\mathbb R\\
 s&\mapsto&  T^{-1}({1\over s})
\end{array}\right.
\end{eqnarray}
which satisfies the following implicit equation 
$\displaystyle 
 T(u(s))=\frac 1s.
$ The hyperbolic metric with constant curvature $-1$ correspond to the function $u(s) = \log s$.
This function $u$ is of interest since it gives precise estimates (up a \emph{bounded} term) of the distance between points lying on the same  horosphere  ${\mathcal H}_t:= \{(x, t): x \in \mathbb R^{N-1}\}$ where $t \in \mathbb R$ is fixed.  
Namely, the  distance between ${\bf x}_t:=(x, t)$
 and ${\bf y}_t:= (y, t)$ for the   metric $ T^2(t){\rm d}x_{ }^2$ induced by $g$ on ${\mathcal H}_t$ is equal to $T(t)\Vert x-y\Vert_{ }$. Hence,   this distance  equals $1$ when $t=u( \Vert x-y\Vert_{ })$ and the union of the 3 segments $[{\bf x}_0, {\bf x}_t],
 [{\bf x}_t, {\bf y}_t] $ and $[{\bf y}_t, {\bf y}_0]$ lies at a bounded distance of the hyperbolic geodesic joigning  ${\bf x}_0$ and ${\bf y}_0$ (see  \cite{DOP}, lemme 4)  : this readily implies that $d({\bf x}_0, {\bf y}_0)-2u( \Vert x-y\Vert)$ is bounded.

The ``curvature''  function $K_g$  may be expressed in term of $u$ as follows:
\begin{equation}\label{curvature}
K_g(u(s)):=-\frac {T''(u(s))}{T(u(s))}=
-\frac {2 u'(s)+s u''(s)}{s^2(u'(s))^3}.
\end{equation}

For any $\alpha \geq 0$,  let us consider the  non decreasing $C^2$-function $ u=u_{\alpha} $ from $ \mathbb R^{*+}$
to
$\mathbb R$ such that 
$$
(i)\quad u_\alpha (s)=\log s \  \mbox{\rm if}
\ 0<s\leq 1\quad {\rm and} \quad (ii)\quad 
u_\alpha (s)= \log s+ \alpha \log \log s   \quad  \mbox{\rm if} \quad
s \geq s_{\alpha }   
$$
for some constant $s_{\alpha } >1$ to be chosen in the following way.
Using formula (\ref{curvature}) and following  Lemma 2.2 in \cite{P},  for any $A\in ]0, 1[$, one may choose $s_{\alpha}  >1$ in such a way 
 the metric $  g_\alpha =T_{u_{\alpha} }^2(t){\rm d}x_{ }^2+{\rm d}t ^2$ on $\mathbb R^{N-1}\times \mathbb R$  has  pinched negative  curvature on $X$,  bounded from above by $-A^2$. 
 Let us emphasize that this metric co\"incides with the hyperbolic one on the subset $\mathbb R^{N-1}\times \mathbb R^-$ and that we can enlarge this subset  shifting the metric $g_\alpha $   along the  axis $\{0\}\times \mathbb R$   as  far as we want  (see \cite{P} $\S$ 2.2).

  Now, let $\mathcal P $ be a discrete  group of isometries  of
$  \mathbb R^{N-1}$
  spanned  by $k$ linearly independent translations $p_{\vec{\tau}_1}, \cdots, p_{\vec{\tau}_k}$ in $  \mathbb R^{N-1}$.  For any ${\bf n}=(n_1, \cdots, n_k) \in \mathbb Z^k$, we set  $\vec{\bf n}= n_1\vec{\tau}_1+ \cdots +n_k\vec{\tau}_k$. The  translations $p_{\vec{\bf n}}$ are also isometries of $(\mathbb R^N, g_\alpha)$  and the  corresponding Poincar\'e series of $\mathcal P$  is given by, up to finitely many terms,
\begin{eqnarray*}
 P_{\mathcal P}(s) = \sum_{\Vert\vec{\bf n}\Vert > s_{\alpha, \beta}}  e^{-sd({\bf o}, p_{\vec{\bf n}}\cdot{\bf o})}&=& \sum_{ \Vert\vec{\bf n}\Vert > s_{\alpha, \beta}} 
 {e^{-2su( \Vert\vec{\bf n}\Vert )-s O( 1)}}\\
 &=&
  \sum_{ \Vert\vec{\bf n}\Vert > s_{\alpha, \beta}} 
{e^{-s O(1)}\over 
 \Vert\vec{\bf n}\Vert^{2s}
\Bigl(\log\Vert \vec{\bf n}\Vert \Bigr)^{2s\alpha} 
}. \end{eqnarray*}
Thus, the  Poincar\'e exponent   of $\mathcal P$  equals   ${k/2}$  and  $\mathcal P$ is convergent if and only if $\alpha >1$.

\begin{remark}
We can construct other  similar metrics as follows. For $\alpha>1$, $\beta>0$, there exists $s_{\alpha, \beta}>0$ and $u_{\alpha,\beta} : (0, +\infty)\to \mathbb R$ such that

(i) $\quad u_{\alpha,\beta}(s)=\log s\quad $ if  
$\quad 0<s\leq 1$,

(ii) $\quad 
u_{\alpha,\beta}(s)= \log s + \alpha \log\log s + \beta \log \log \log s \quad $ if  
$\quad s\geq s_{\alpha,\beta}$,

(iii) $K_g(u(s))\leq -A$. 
\\
Hence, the Poincar\'e series of the parabolic subgroup $\mathcal P$  with respect to the metric $ g_{\alpha, \beta} = T^2_{u_{\alpha,\beta}}(t)^2 {\rm d}x^2 + {\rm d}t^2$ is  given by, up to finitely many terms,

\begin{eqnarray*}
 P_{\mathcal P}(s) = \sum_{ \Vert\vec{\bf n}\Vert > s_{\alpha, \beta}}  e^{-sd({\bf o}, p_{\vec{\bf n}}\cdot{\bf o})}&=& \sum_{ \Vert\vec{\bf n}\Vert > s_{\alpha, \beta}} 
 {e^{-2su( \Vert\vec{\bf n}\Vert  )-s O(1)}}\\
 &=&
  \sum_{ \Vert\vec{\bf n}\Vert > s_{\alpha, \beta}} 
{e^{-s O(1)}\over 
 \Vert\vec{\bf n}\Vert^{2s}
\Bigl(\log \Vert\vec{\bf n}\Vert\Bigr)^{2s\alpha} \Bigl(\log\log\Vert\vec{\bf n}\Vert\Bigr)^{2s\beta}
}. \end{eqnarray*}
\end{remark}

\medskip
This implies that $\mathcal P$ converges as soon as $\alpha>1$ but it  is not enough to ensure that $\mathcal P$   satisfy   hypothesis (\ref{influent}). In the next paragraph, we   present new metrics $g_\alpha$, close to those presented in the present section, for which it holds
$$d(o, p_{\vec{\bf n}}\cdot  o) = 2\left (\log  \Vert\vec{\bf n}\Vert  + \alpha \log \log  \Vert\vec{\bf n}\Vert \right) + C + \epsilon(n),$$
where $C\in \mathbb R$ is a constant   and $\displaystyle \lim_{n\to +\infty} \epsilon(n) = 0$. 

 \subsection{ On convergent parabolic groups satisfying condition (\ref{influent})}

Let us fix $N=2, \alpha>1$ and  a slowly varying function $L: [0, +\infty[\to \mathbb R^{*+}$. We   construct in this section a metric $g=g_{\alpha, L}=T^2(t){\rm d}x^2 + {\rm d} t^2$ on $\mathbb R\times \mathbb R$ such that the group spanned  by the translation $(x,t) \mapsto (x+1, t)$ satisfies our hypothesis (\ref{influent}). The generalization to higher dimension is immediate.

 For any  real $t$ greater than some ${\mathfrak a}>0$ to  be chosen, let us set
$$T(t)=T_{  \alpha, L}(t) = e^{-t}{t^\alpha\over L(t)}.$$
 Without loss of generalities, we assume that $L$ is $C^{\infty}$ on $\mathbb R^+$ and  its derivates $L^{(k)}, k \geq 1$, satisfy  
$
L^{(k)}(t)\longrightarrow 0
$
and 
$
\displaystyle \frac{L'(t)}{L(t)} \to 0
$ as $t\to +\infty$ (\cite{BGT}, Section 1.3). Furthermore, for any $\theta>0$, there exist $t_\theta \geq 0$ and $C_\theta \geq 1$ such that  for any $t\geq t_\theta$
\begin{equation} \label{majorationslowlyvarying}
 {1\over C_\theta t^\theta} \leq L(t)\leq C_\theta t^\theta.
 \end{equation}
Notice that 
$\displaystyle 
  -\frac{T ''(t)}{T (t)} = -\left(1-{2\alpha\over t}+L'(t)\right)^2+\left({\alpha \over t^2}+L''(t)\right)   <0 $ for $t \geq{\mathfrak a}$.

We assume that $0<A<1<B$ and, following Lemma 2.2 in \cite{P}, extend  $T_{\alpha, L}$ on $\mathbb R$ as follows. 
\begin{lemma}\label{talpha}
There exists  $\mathfrak a= \mathfrak a_{\alpha, L}>0$  such that the map $T=T_{\alpha, L} : \mathbb R \to (0, +\infty)$ defined by
\begin{itemize}
\item $T (t) = e^{-t}\quad $ for  $\quad t\leq 0$, 
\item $T(t)= e^{-t} {t^\alpha\over L(t)}\quad $ for  $\quad  t\geq \mathfrak a_{\alpha, L}$,
\end{itemize}
admits a decreasing and 2-times continuously differentiable  extension on $\mathbb R$ satisfying the following inequalities
$$
-B \leq K(t) = -\frac{T''(t)}{T(t)}\leq -A<0.
$$
\end{lemma} 
Notice that  this property holds for  any $t'\geq t_{\alpha, L}$;  
 {\bf for technical reasons} (see Lemma \ref{groupeconvergentHypotheses}),  we  assume   that $\mathfrak a > 4\alpha$.
A direct computation  yields the following estimate for the   function $u=u_{\alpha, L}$  given by the implicit equation  (\ref{eq:FonctionU}).
\begin{lemma}\label{lem:UAsymp}
Let $u=u_{\alpha, L} : (0, +\infty) \to \mathbb R$ be such  $\displaystyle T_{\alpha, L}(u(s)) = \frac 1 s$  for any $s>0$. Then 
$$u(s)= \log s +\alpha\log \log s-\log L(\log s)+\epsilon(s)$$
 with $\epsilon(s) \to 0$ as $s \to +\infty$.
\end{lemma}

We now consider  the group $\mathcal P$ spanned  by  the translation $p$ of vector $\vec{i} = (1,0)$ in $\mathbb R^2$; 
 the map $p$ is an  isometry of $(\mathbb R^2, g_{\alpha, L})$ which  fixes the point $x=\infty $. By Lemma \ref{lem:UAsymp}, it holds  
 $$
 d({\bf o}, p^n\cdot{\bf o})=2\Bigl(\log n +\alpha\log \log n-\log L(\log n)\Bigr)
 $$
  up to a bounded term. Hence, the   group    $\mathcal P$ has critical exponent $\frac 1 2$; furthermore, it  is convergent since $\alpha>1$. $^($\footnote{Notice that the group $\mathcal  P$ also converges when $\alpha = 1$  and 
  $\displaystyle\sum_{n \geq 1}{L(n)\over n} <+\infty;$ this situation is not  explore here.}$^)$
  The following proposition   ensures that $\mathcal P$ satisfies  hypothesis (\ref{influent}); in other words, the ``bounded term'' mentioned above tends to $0$ as $n \to +\infty$.

 \begin{prop} \label{prop:PAsymp} 
The parabolic  group $\mathcal P
=\langle p \rangle$  on $(\mathbb R^2, g_{\alpha, L})$ satisfies the following property: for any $n\in \mathbb N$,
$$
d({\bf o}, p^n\cdot{\bf o})=  2\Bigl(\log n +\ \alpha\log \log n -\log L(\log n)\Bigr)  + \epsilon(n)
$$
with $\displaystyle \lim_{n \to +\infty} \epsilon(n)=0$.  In particular, if $\alpha >1$, then $\mathcal P$   is convergent with respect to $g_\alpha$.
 \end{prop}

Let $\mathcal H = \mathbb R \times [0, +\infty)$ be the upper half plane $\{(x, t) \mid  t\geq 0\}$  and $\mathcal H/\mathcal P$ the quotient cylinder endowed with the metric $g  _{\alpha, L} =T  _{\alpha, L}(t)^2{\rm d}x^2+{\rm d}t^2$.  We   do not estimate directly the distances $d({\bf o}, p^n\cdot{\bf o})$, since the metric $g  _{\alpha, L}$ is not known explicitely for $t\in [0, \mathfrak a]$. Let us introduce the point ${\bf a} = (  0, {\mathfrak a})\in \mathbb R^2$. The union of the three geodesic segments $[{\bf o}, {\bf a}], [{\bf a}, p^n\cdot{\bf a})]$ and $[p^n\cdot {\bf a}, p^n\cdot {\bf o}]$ is a quasi-geodesic; more precisely, since $d({\bf o}, {\bf a})= d(p^n\cdot {\bf o}, p^n \cdot {\bf a})$ is fixed  and $d({\bf a}, p^n\cdot{\bf a})\to +\infty$,   the following statement holds.
\begin{lemma}
Under the previous notations, 
$$\lim_{n\to +\infty} d({\bf o}, p^n \cdot{\bf o})  - d({\bf a}, p^n\cdot{\bf a}) = 2\mathfrak a.$$
\end{lemma}
Proposition \ref{prop:PAsymp}   follows from the following lemma.

 \begin{lemma}
\label{groupeconvergentHypotheses} Assume that $\mathfrak a \geq  4 \alpha$. Then
$$
d({\bf a}, p^n\cdot{\bf a})=  2(\log n +\ \alpha\log \log n -\log L(\log n)- \mathfrak a )+ \epsilon(n)
$$
with $\displaystyle \lim_{n \to +\infty} \epsilon(n)=0$. 

 \end{lemma} 
Proof. 
%
Throughout this proof, we work on the  upper half-plane $\mathbb R \times [\mathfrak a, +\infty[$ whose points are denoted $(x, \mathfrak a+t)$ with $x \in \mathbb R$ and $ t\geq 0$; we set  
  $${\mathcal T}(t) = T  _{\alpha}(t+\mathfrak a)=  e^{-\mathfrak a-t}{(t+\mathfrak a)^\alpha\over L(t+\mathfrak a)}.$$
In these   coordinates, the quotient cylinder $\mathbb R \times [\mathfrak a, +\infty[/\mathcal P$   is a surface of revolution endowed with the metric  ${\mathcal T}(t)^2{\rm d}x^2+{\rm d}t^2$. 
For any $n \in \mathbb Z$, denote $h_n$ the maximal height at which the geodesic segment $\sigma_n=[{\bf a}, p^n\cdot {\bf a}]$  penetrates inside the upper half-plane $\mathbb R \times [\mathfrak a, +\infty[$; it tends to $+\infty$ as $n \to \pm  \infty$.
The relation between 
 $n, h_n$ and    $d_n:=d({\bf a}, p^n\cdot{\bf a})$ may be deduced from 
the
Clairaut's relation (\cite{DC}, section 4.4, Example 5) :
$$ 
   {n\over 2}={\mathcal T}(h_n)\int_{0}^{\mathfrak  h_n}{{\rm d}t\over 
{\mathcal T}(t)\sqrt{{\mathcal T}^2(t)-{\mathcal T}^2(h_n)}}
\quad {\rm and} \quad   d_n=2\int_{0}^{  h_n}{{\mathcal T}(t){\rm d}t\over 
\sqrt{{\mathcal T}^2(t)-{\mathcal T}^2(h_n)}}.
$$
These identities may be rewritten as
$$ 
  {n\over 2} ={1\over {\mathcal T}(h_n)}\int_{0}^{\mathfrak  h_n}{f_n^2(s){\rm d}s\over 
 \sqrt{1-f_n^2(s)}}
\qquad {\rm and} \qquad    \quad d_n=2 h_n+ 2 \int_{0}^{ h_n}
\Bigl({1\over 
\sqrt{1-f_n^2(s)}}-1\Bigr) {\rm d}s
$$
where $\displaystyle f_n(s):= {{\mathcal T}(h_n)\over {\mathcal T}(h_n-s)}1_{[0, h_n]}(s).$

First,  for any $s \geq 0$, the quantity $\displaystyle  
{f_n^2(s)\over 
 \sqrt{1-f_n^2(s)}}$ converges towards $\displaystyle {e^{-2s}\over\sqrt{1-e^{-2s}}}
 $ as $n \to +\infty$.  In order to use the dominated convergence theorem, we  need the following property.

\begin{fact}\label{majorationf}
There exists $n_0>0$ such that for any $n\geq n_0$ and any $s\geq 0$,
$$
0\leq f_n(s)\leq f(s):= 
e^{-s/2} 
$$
\end{fact}

\noindent Proof. 
Assume first  $h_n/2\leq s\leq h_n$; taking $\theta=\alpha/2$ in (\ref{majorationslowlyvarying}) yields 
\begin{eqnarray*}
0\leq f_n(s)&=& 
\left(\mathfrak a + h_n\over \mathfrak a + h_n - s\right)^\alpha  { L( \mathfrak a + h_n - s)\over L( \mathfrak a + h_n)} e^{-s}
\\
&\leq&
C_{\alpha/2}^2{(\mathfrak a + h_n)^{3\alpha/2}\over (\mathfrak a + h_n - s)^{\alpha/2} } e^{-s}
\\
&\leq&
{C_{\alpha/2}^2  \over
\mathfrak a^{\alpha/2}}
(\mathfrak a + h_n)^{3\alpha/2}
 e^{-s}
 \\
&\leq&
{C_{\alpha/2}^2  \over
\mathfrak a^{\alpha/2}}
(\mathfrak a + h_n)^{3\alpha/2}e^{-{h_n\over 4}} e^{-{s\over 2}} \leq e^{-{s\over 2}}
\end{eqnarray*}
where the last inequality holds if $h_n$ is great enough, only depending on $\mathfrak a$ and $\alpha$.

Assume now $0\leq s \leq h_n/2$; it holds  $\displaystyle {1\over 2}\leq {\mathfrak a+h_n-s\over \mathfrak a +h_n}\leq 1 $ and 
  $0\leq {s\over \mathfrak a +h_n}\leq \min({1\over 2}, {s \over \mathfrak a})$. 
  Recall that $L'(t)/L(t)\to 0$ as $t\to +\infty$ and $0\leq {1\over 1-v }\leq e^{2v}$ for $0\leq v \leq {1\over 2}$; hence, for any 
  $\varepsilon >0$ and $n$ great enough (say $n \geq n_\varepsilon$), there exists $s_n\in (0,s)$ such that
 \begin{eqnarray*}
0\leq f_n(s)&=&  { L( \mathfrak a + h_n - s)\over L( \mathfrak a + h_n)} 
\left(\mathfrak 1\over \mathfrak 1 - {s\over \mathfrak a + h_n}\right)^\alpha e^{-s} 
\\
&\leq& \left(1 - s\frac{L'(a+h_n - s_n)}{L(a+h_n)}\right) e^{-(1-{2\alpha \over \mathfrak a})s} \\
& \leq & (1 + \epsilon s)e^{-(1-{2\alpha \over \mathfrak a})s}\\
&\leq&  e^{-(1-\varepsilon -{2\alpha \over \mathfrak a})s}.
 \end{eqnarray*}
 Consequently, fixing $\epsilon >0$ in such a way   $\displaystyle 2{\alpha\over\mathfrak a} + \epsilon \leq \frac 1 2$, it yields $0\leq f_n(s)\leq e^{-s/2}$ for $n$ great enough. 
\rightline{$\Box$}
Therefore, 
$$
0\leq {f_n^2(s)\over 
 \sqrt{1-f_n^2(s)}}
 \leq  F(s):= {f^2(s)\over 
 \sqrt{1-f^2(s)}}
$$
where  the function $F$ is integrable on $\mathbb R^+.$
By the dominated convergence theorem, it yields 
$${n\over 2}={1+\epsilon(n)\over {\mathcal T}(h_n)} \int_0^{+\infty}  {e^{-2s}\over\sqrt{1-e^{-2s}}}{\rm d} s
={1+\epsilon(n)\over  \mathcal {\mathcal T}(h_n)}.$$
Consequently   $h_n= \log n +\alpha \log \log n -\log L(\log n)-\log 2 -\mathfrak a +\epsilon(n)$.

Similarly $\displaystyle \lim_{n\to +\infty}
\int_0^{h_n}
\Bigl({1\over 
\sqrt{1-f_n^2(s)}}-1\Bigr) {\rm d}s=
\int_0^{+\infty}
\Bigl({1\over 
\sqrt{1-e^{-2s}}}-1\Bigr) {\rm d}s=  \log 2, 
$
which yields   $$d_n= 2(\log n +\alpha \log \log n -\log L(\log n) -\mathfrak a) +\epsilon(n).$$

\rightline{$\Box$}

  The Poincar\'e exponent of $\mathcal P$  equals $1/2 $ and, as $R\to +\infty$, 
$$
\sharp\{p\in \mathcal P\mid 0\leq d({\bf o}, p\cdot {\bf o})<R \} \sim   e^{ R/2}{ L(R)\over (R/2)^{\alpha} }.
$$
Hence, for any $\Delta >0$, 
$$
\sharp\{p\in \mathcal P\mid R\leq d({\bf o}, p\cdot {\bf o})<R+\Delta\} \sim {1\over 2}\int_{R}^{R+\Delta} e^{t/2}{L(t)\over (t/2)^{\alpha}}  {\rm d}t\quad{\rm as}\quad R\to+\infty
$$
and
$$
 \lim_{R\to +\infty} { R ^{\alpha}  \over L(R)}
 \sum _{ \stackrel{p \in \mathcal P}{R\leq d({\bf o}, p\cdot {\bf o})<R+\Delta}} 
 e^{-{1\over 2} d({\bf o}, p\cdot {\bf o})}  =2^{\alpha -1}\Delta
$$
which is precisely  Hypothesis  {\ref{influent}.
%

 \subsection{On the existence of  non elementary exotic groups} 
 
Explicit constructions of exotic groups, i.e. non-elementary groups $\Gamma$ containing a parabolic $\mathcal P$ whose Poinacr\'e exponent equals $\delta_\Gamma$, have been detailed in several papers; first in \cite{DOP}, then in  \cite{P}, \cite{DPPS} and \cite {V}. Let us describe them in the context of the metrics $g=g_{\alpha; L}$ presented above.

 For any $a>0$ and $t\in \mathbb R$, we   write
$$
T_{\alpha, L, a} = \left\{\begin{array}{ccc}
					e^{-t} & \mbox{if} & t\leq a\\
					e^{-a}T_{\alpha, L}(t-a) & \mbox{if} & t\geq a
					\end{array}\right.,
$$
where $T_{\alpha, L}$ is defined in the previous paragraph. As in \cite{P}, we consider the metric on $\mathbb R^2$ given   by $\displaystyle g_{\alpha, L, a} = T_{\alpha, L, a}^2(t){\rm d}x^2 + {\rm d}t^2$. It is a complete smooth metric, with pinched negative curvature, and which equals the hyperbolic one  on $\mathbb R\times (-\infty, a)$. Note that $g_{\alpha, L,  0} = g_{\alpha, L}$ and $g_{\alpha, L, +\infty}$ is the hyperbolic metric on $\mathbb H^2$. Note the  previous subsection, for any $a\in (0, +\infty)$ and any $\tau\in \mathbb R^*$,  a  parabolic group of the form $\mathcal P = <(x,t) \mapsto (x+\tau, r)>$ is convergent. This allows   to reproduce the construction of a non-elementary group given in \cite{DOP} and \cite{P}.

Let $h$ be a hyperbolic isometry of $\mathbb H^2$ and $p$ be a parabolic isometry in Schottky position with $h$ (cf next section for a precise definition). They generate a free group $\Gamma = <h,p>$ which acts discretely without fixed point on $\mathbb H^2$.  Up to a global conjugacy, we can suppose that $p$ is $(x,t) \mapsto (x+\tau, t)$ for some $\tau\in \mathbb R^*$. The surface $S = \mathbb H^2/\Gamma$ has a cusp, isometric to $\mathbb R/ \tau\mathbb Z\times (a_0, +\infty)$ for some $a_0>0$. Therefore, we can  replace in the cusp the hyperbolic metric by $g_{\alpha, L, a}$ for any $a\geq a_0$; we also denote $ g_{\alpha, L, a}$ the lift of $g_{\alpha, L,  a}$ to $\mathbb R^2$. 

For any $n\in \mathbb Z^*$, the group $\Gamma_n = <h^n, p>$ acts freely by isometries on $(\mathbb R^2, g_{\alpha, L, a})$. It is shown in   \cite{DOP}  that,   for $n>0$ great enough, the group $\Gamma_n$  also converges. This provides a family of examples for Theorem \ref{theo:Comptage}. By \cite{P}, if $\Gamma_n$ is convergent for some $a_0>0$, then there exists $a^*>a_0$ such that for any $a\in [a_0, a^*)$, the group $\Gamma_n$ acting on $(\mathbb R^2, g_{\alpha, L,  a})$ is   convergent, whereas for $a>a^*$,  it has finite Bowen-Margulis measure and hence diverges. In some sense, the case $a = a^*$ is ``critical'';  it is proved  in \cite{P} that $\Gamma$  also diverges in this case. With additive hypotheses on the   tail   of the Poincar\'e series associated to the factors $\Gamma_j, 1\leq j \leq p$ of $\Gamma$, P. Vidotto  has obtained a precise estimate of the orbital function of $\Gamma$ in the case when its Bowen-Margulis measure is infinite \cite {V} ; this is the analogous of  Theorem \ref{theo:Comptage},  under slightly more general assumptions.
 
 In  \cite{DPPS}, the authors propose another approach   based on a ``strong'' perturbation of the metric inside the cusp.  Starting from a $N$-dimensional  finite volume hyperbolic manifold with cuspidal ends, they modify the metric far inside one end  in  such a way the corresponding parabolic group is convergent with  Poincar\'e exponent $>1$ and  turns the fundamental group of the manifold into a convergent group;  in this construction,  the   sectional curvature of the new metric along certain planes is $<-4$ far inside the modified cusp.


\section{Schottky groups: generalities and coding}\label{sec:Schottky}
From now on, we fix two integers $p\geq 1$ and $q\geq 0$ such that $\ell:=p+q\geq 2$ and consider a Schottky group $\Gamma$ generated by $\ell $  elementary  groups 
 $\Gamma_1, \ldots, \Gamma_\ell$ of  isometries of $X$.   These  elementary groups are   in Schottky position, i.e.   there
exist disjoint closed
sets
$F_j$ in  $\partial X$ such that, for any $1\leq j\leq \ell$
 $$
\Gamma _j^*(\partial X\setminus F_j) \subset F_j.
$$
The group $\Gamma$ spanned  by the   $\Gamma_j, 1\leq j \leq \ell,$ is called the
Schottky product of the ${\Gamma_j}$'s  and denoted $\Gamma=\Gamma_1\star 
\Gamma_2\star \cdots \star
\Gamma_\ell$. 

In this  section, we present general properties of $\Gamma$. In particular, we do not require that conditions {\bf H1}, {\bf H2} and {\bf H3} hold; these hypotheses are only  needed in the   
last section  of this paper.

By
the Klein's tennis table criteria,  $\Gamma$ is the free product of
the groups $\Gamma _i$; any element in $\Gamma$ can be   uniquely written 
as the product
$$ \gamma = a_1\dots a _k$$ for some  $a_j \in \cup\Gamma_j^*$ with the 
property  that  no two consecutive elements $a_j$ belong to the same 
group. The set $\mathcal A= \cup\Gamma_j^*$ is  called the {\it alphabet } of 
$\Gamma$, and $ a_1, \dots, a_k$ the {\it letters} of $\gamma$. The number 
$k$ of letters  is  the {\it
symbolic length }
of
$\gamma$; let us denote $\Gamma(k)$ the set of elements of $\Gamma$ with symbolic length $k$. The last
letter of $\gamma$   plays a special role, and the index of the group it 
belongs to   be   denoted  by $l_\gamma$.  Applying Fact 2 one 
gets\

\begin{property} \label{triangle} There exists a constant $C>0$ such that
 $$d({\bf o}, \gamma.{\bf o})-C\leq B_x(\gamma^{-1}.{\bf o}, {\bf o})\leq d({\bf o}, \gamma.{\bf o})$$
 for 
any $\gamma\in \Gamma= \star _i\;\Gamma _i$ and any  $x\in 
\cup_{i\not=l_\gamma}F_i$.
\end{property}
\noindent 
This fact implies in particular   the   following crucial 
contraction
property \cite{BP}.
\begin{prop}\label{contraction} There exist a real number $r\in ]0, 1[$ 
and
$C>0$ such that for any
$\gamma $ with symbolic length $n\geq 1$ and any $x$ belonging to the closed set
$\cup_{i\not= i(\gamma
)}F_i$ one has $$\vert\gamma'(x)\vert
\leq Cr^n.$$
\end{prop}
%
%
%

The following statement, proved in \cite{BP},    provides a coding of the limit set $\Lambda(\Gamma)$ but   the $\Gamma$-orbits of the fixed points of the generators.
 
\begin{prop}\label{codagelimitset}  Denote by $\Sigma^+$ the set of 
sequences $(a
_n)_{n\geq 1}$
for which each letter  $a _n$ belongs  to the alphabet $\mathcal A= \cup\Gamma_i
^*$ and such  that
no two consecutive letters belong to the same group (these sequences are 
called admissible). Fix a point $x_0$ in
$\partial X\setminus F$. Then
\begin{enumerate}
\item[(a)] For any
${\bf a}= (a_n)_{n \geq 1}\in
\Sigma^+$, the sequence
$(a_1\dots  a_n\cdot x_0)_{n \geq 1}$ converges to a point $\pi ({\bf a})$ in the limit set of
$\Gamma$, independent on the choice of $x_0$. 
\item [(b)] The map
$\pi : \Sigma  ^+\to \Lambda (\Gamma )$  is
one-to-one and  $\pi(\Sigma^{+})$ is contained in the radial limit 
set of $\Gamma$.
\item [(c)] The complement of $\pi(\Sigma  ^+)$ in the limit set of  
$\Gamma $
equals the
$\Gamma$-orbit of the union of the  limits sets $\Lambda (\Gamma _i)$
\end{enumerate}
\end{prop}

From now on, we consider a Schottky product  group $\Gamma$.  Thus, up to a denumerable set of  points, the limit set of  
$\Gamma$ coincides with
$ \pi (\Sigma^+)$. For any $1\leq i \leq \ell$, let $\Lambda_i = \Lambda\cap F_i$
be the closure of the set  of those limit points with first letter in 
$\Gamma_i $ (not to be confused with the limit set of $\Gamma_i$). The 
following
description of $\Lambda=\Lambda(\Gamma)$  be useful:

a) $\Lambda$ is  the finite union of the   sets $\Lambda_i$, 

b) the closes sets $\Lambda _i, 1\leq i \leq \ell,$ are pairwise disjoints,

c) each of these sets is partitioned into a countable number of
subsets with disjoint closures\ :
$$
\Lambda_i= \cup_{a \in \Gamma^*_i}\cup_{j\not= i}\;
a.\Lambda_j\ .
$$

Now, we enlarge the set $\Lambda$ in order to take into account the finite  admissible  words. We fix a point $x_0\notin \cup_jF_j$. There exists a one-to-one correspondence between $\Gamma\cdot x_0$ and $\Gamma$; furthermore, the point $\gamma\cdot x_0\in F_j$  for any $\gamma \in \Gamma^*$  with first letter in $\Gamma_j$. We  set $\widetilde \Sigma_+= \Sigma^+\cup \Gamma$ and  notice that, by the previous Proposition, 
the natural map $\pi: \widetilde \Sigma_+ \to \Lambda(\Gamma)\cup \Gamma\cdot {x_0}$ is one-to-one  with image $\pi(\Sigma^+)\cup \Gamma \cdot x_0$. Thus we introduce the following notations

a) $\tilde \Lambda= \Lambda \cup \Gamma \cdot x_0;$

b) $\tilde \Lambda_i= \tilde \Lambda\cap F_i$ for any $1\leq i \leq \ell$.

\noindent The set $\tilde \Lambda$ is the disjoint union of $\{x_0\}$ and the sets $\tilde \Lambda_i, 1\leq i \leq \ell$; furthermore, each $\tilde\Lambda_i$ is partitioned into a countable number of
subsets with disjoint closures:
$$
\tilde \Lambda_i= \cup_{a \in \Gamma^*_i}\cup_{j\not= i}\;
a\cdot \tilde \Lambda_j\ .
$$

The cocycle $b$  defined in  (\ref{cocycleb})   play a central role in the sequel.  In order to calculate the distance between two points of the orbit $\Gamma \cdot {\bf o}$, we consider an extension $\tilde b$ of this  cocycle defined as follow on $\tilde \Lambda$: for any $\gamma \in \Gamma$ and $x \in \tilde \Lambda$,
$$
\tilde b(\gamma, x):=
 \Bigl\{
 \begin{array}{lllll}
  b(\gamma, x)= \mathcal B_x(\gamma^{-1} {\bf o}, {\bf o})& {\rm if} &  x\in 
 \Lambda;&\ &  \\
d(\gamma^{-1} \cdot {\bf o}, g\cdot {\bf o})-d({\bf o}, g\cdot{\bf o})
& {\rm  if} &  x=g\cdot x_0& {\rm for \ some} & g \in \Gamma.
\end{array}
\Bigr.
$$
The cocycle equality (\ref{cocycleb}) is still valid for the function $\tilde b$; 
furthermore, if $\gamma\in \Gamma$ decomposes as $\gamma =a_1\cdots a_k$, then
$$
d({\bf o}, \gamma \cdot {\bf o}) =  b(a_1, \gamma_2\cdot x_0)  +b(a_2, \gamma_3\cdot  x_0)+\cdots+ b(a_k, x_0),  
$$
where  $\gamma_l= a_l\cdots a_k$ for $2\leq l\leq k$.

\section{ On the Ruelle operators $\mathcal L_ {s}, s \in \mathbb R$}\label{sec:Ruelle}
In this  
section, we  describe the main properties   of 
the  transfer operators   $\mathcal L_ {s}, s \in \mathbb R,$  defined
formally 
 by:  for any function $\phi:  
 \tilde   \Lambda \  \to 
\mathbb C$ and $x \in  \tilde \Lambda$,
$$\mathcal L_ {s}\phi(x)= 
\sum_{\gamma\in \Gamma(1)} {\bf 1}_{x \notin  \tilde   \Lambda_{l_\gamma}} e^{-s\tilde b(\gamma, x)}\phi(\gamma\cdot  x)=\sum_{j=1}^{\ell}
 \sum_{\gamma  \in \Gamma^*_j }  
{\bf 1}_{x \notin {\tilde \Lambda_j}}  e^{-s\tilde b(\gamma, x)}\phi(\gamma\cdot  x).
$$
For any $1\leq j\leq \ell$, the sequence $(\gamma \cdot {\bf o})_{\gamma \in \Gamma_j}$ accumulates on the  fixed point(s) of $\Gamma_j$. So for any $x \notin  \tilde \Lambda_j$, the  sequence $\left(
\tilde b(\gamma, x)- d ({\bf o}, \gamma.{\bf o}
)\right)_{\gamma \in \Gamma_j}$ is bounded uniformly in $x \notin \tilde \Lambda_j$. Therefore the  quantity $\mathcal L_ {s}1(x)$ is well defined as soon as $s\geq \delta:=  \max\{\delta_{\Gamma_j}\mid 1\leq j \leq \ell\}$.
The powers of $\mathcal L_ {s}, s\geq \delta,$ are formally given by: for any $k \geq 1$,  any function $\phi:  
\tilde \Lambda \  \to 
\mathbb C$ and any $x \in \tilde \Lambda$, 
$$
\mathcal L^k_ {s}\phi(x)= 
 \sum_{\gamma  \in \Gamma(k)}
{\bf 1}_{x \notin {\tilde \Lambda_j}}  e^{-s\tilde b(\gamma, x)}\phi(\gamma\cdot  x).
$$
It  is easy to check that the operator    $\mathcal L_ {s}, s \geq \delta$, act   on $(C(\tilde \Lambda), \vert \cdot \vert _\infty)$; we denote 
$  \rho_s(\infty)$
it spectral radius on this space.
%

\subsection{ Poincar\'e series versus  Ruelle operators}

  By the ``ping-pong dynamic'' between the subgroups $\Gamma_j, 1\leq j\leq \ell$,  and Property  \ref{triangle},   we easily check  that the difference
 $
 \tilde b(\gamma, x)-d({\bf o}, \gamma\cdot{\bf o})
 $
 is bounded uniformly in $k\geq 0,  \gamma \in \Gamma(k)$ and $x \notin \tilde \Lambda_{l_\gamma}$. Consequently, there exists  a constant $C >0$ such that, for any $x \in \tilde \Lambda$,  any $k \geq 1$ and any $s \geq \delta$,
 $$
\mathcal L^k_ {s}1(x)  \stackrel{c}{\asymp} \sum_{\gamma \in \Gamma(k)} e^{-s d({\bf o}, \gamma \cdot{\bf o})}
 $$
 where $A\stackrel{c}{\asymp}B$ means ${A\over c}\leq B \leq cA$.
Hence, 
\begin{equation}\label{divergence-convergence}
 \displaystyle P _\Gamma(s):= \sum_{ \gamma \in \Gamma} e^{-sd({\bf o}, \gamma\cdot{\bf o})}=+\infty 
 \quad \Longleftrightarrow\quad 
 \sum_{k\geq 0} \mathcal L^k_ {s}1(x)=+\infty.
 \end{equation}
In particular
\begin{equation}\label{criticalexponent}
\delta_\Gamma=\sup\{s\geq \delta\mid \rho_s(\infty)\geq 1\}=  \inf\{s\geq \delta\mid \rho_s(\infty)\leq 1\}.
 \end{equation}

It is proved in the next paragraph that $\Gamma$ is convergent if and only if $\rho_\delta(\infty)<1$.

\subsection{ On the spectrum of the operators  $\mathcal L_ {s}, s\geq \delta $}

In order to control the spectral radius (and the spectrum) of the transfer operators $\mathcal L_{s}$, we   study their restriction  
to the space
${\bf Lip}( \tilde \Lambda) $
of Lipschitz functions from $ \tilde \Lambda$ to $\mathbb C$
defined by
$${\bf Lip}( \tilde \Lambda)=\{\phi \in C( \tilde \Lambda);\; \Vert\phi\Vert =
|\phi|_{\infty}+[\phi] <+\infty\}$$
where $\displaystyle [\phi] =\sup_{0\leq i\leq p}
\sup_{\stackrel{x, y \in  \tilde \Lambda_j}{x\neq y}}{|\phi(x)-\phi(y)|\over D(x, y)  }$ is  the Lipschitz coefficient 
of $\phi$   on $(\partial X, D)$. 

The space  $({\bf Lip}( \tilde \Lambda),\Vert.\Vert )$ is a  Banach space and 
the identity
map from
$({\bf Lip}( \tilde \Lambda), \Vert.\Vert )$
into $(C( \tilde \Lambda), |.|_{\infty})$ is
 compact.
It is proved in \cite{BP} that  the operators   $\mathcal L_ {s}, s \geq \delta$, act both  on $(C(   \Lambda), \vert \cdot \vert _\infty)$ and   $({\bf Lip}( \Lambda), \Vert \cdot \Vert    )$; P. Vidotto has extended in \cite{V} this property to the Banach spaces $(C(   \tilde \Lambda), \vert \cdot \vert _\infty)$ and   $({\bf Lip}(\tilde  \Lambda), \Vert \cdot \Vert    )$. We denote $ \rho_s$ the spectral radius of $\mathcal L_ {s}$ on    ${\bf Lip}( \tilde \Lambda)$; in the following proposition, we state the  spectral properties of the $\mathcal L_s$  we need in the present paper.
\begin{prop}\label{resumeDCDS} We assume $\ell=p+q\geq 3\ ^($\footnote{Recall that $\ell\geq 2$ since $\Gamma$ is non-elementary. When $\ell = 2$, the real $-\rho_s$ is also a a simple eigenvalue of $\mathcal L_s$; a similar statement to Proposition \ref{resumeDCDS} holds  for the restriction of $\mathcal L_s$ to each space ${\bf Lip}( \tilde \Lambda_i), i=1, 2$  \cite{BP}.
}$^)$. For any $s \geq \delta$,
\begin{enumerate}
\item $ \rho_s=\rho_s(\infty);$
\item  $ \rho_s$ is  a simple eigenvalue of    $\mathcal L_s $ acting on ${\bf Lip}( \tilde \Lambda)$ and the associated eigenfunction $h_s $ is non negative on $ \tilde \Lambda$;
\item there exists $0\leq r <1$ such that the  rest of the spectrum of  $\mathcal L_s $  on ${\bf Lip}( \tilde \Lambda)$  is included in a disc of radius $\leq r \rho_s$.
\end{enumerate}
\end{prop}
\noindent Sketch of the proof.  We refer to \cite{BP} and \cite{V} for the details. 
For any $s\geq 0$ and $\gamma$ in $\Gamma^*$,  let  $w_s( \gamma, .)$ 
be the {\it weight function} defined on $\Lambda(\Gamma)$ by: for any $s \geq \delta$ and $\gamma \in \Gamma$
$$
 w_s(\gamma, x):=
 \Bigl\{
 \begin{array}{lll}
   e^{-s \tilde b(\gamma, x)}& {\rm if} &  x\in 
\tilde \Lambda_j, j\neq l_\gamma,\\
0& {\rm  otherwise.} &  
\end{array}
\Bigr.
$$
Observe that these 
functions satisfy the following  cocycle relation : if $ \gamma_1, \gamma_2 \in \mathcal A$ 
do not belong to the same group $\Gamma_j$, then 
$$
   w_s(\gamma_1\gamma_2, x)=    w_s(\gamma_1, \gamma_2\cdot x) 
   w_s(\gamma_2, x).
$$
Due to this cocycle property, we may write, for any $k \geq 1$, any bounded function $\varphi:  \tilde \Lambda \to \mathbb R$ and any $x \in  \tilde \Lambda$
 $$
\mathcal L_s^k\varphi(x)=\sum_{\gamma \in \Gamma(k)}
   w_s(\gamma, x) \varphi(\gamma\cdot x).
$$ 
In  \cite{BP},  it is proved  that the restriction of the functions
  $   w_s(\gamma, .)$  to the set  $\Lambda$ belong 
 to
 ${\bf Lip}( \Lambda)$ and that  for any $s\geq \delta$ there exists   $C= C(s)>0$ 
such that, for any $\gamma$ in $\Gamma^*$ 
$$
\Vert    w_s(\gamma, 
.)\Vert \leq Ce^{-s d({\bf o}, \gamma.{\bf o})}.
$$
In \cite{V}, Proposition 8.5, P. Vidotto  has proved that the same inequality  holds for  the functions
  $   w_s(\gamma, .)$ on $\tilde \Lambda$. Thus, the operator $\mathcal L_s$ is   bounded on
 ${\bf Lip}( \tilde \Lambda)$ when
$s\geq \delta $.   

In order to describe its spectrum   on ${\bf Lip}( \tilde \Lambda)$, we first write a ``contraction property''  for the iterated operators 
$\mathcal L_s^k$; indeed, 
$$\vert \mathcal L^k_s\varphi(x)-\mathcal L^k_s\varphi(y)\vert
 \leq \sum_{\gamma\in \Gamma(k)} \vert w_s(\gamma,x)\vert
\; \vert\varphi(\gamma\cdot x)-\varphi(\gamma\cdot y)\vert + \sum_{\gamma\in 
\Gamma(k)}
[w_s(\gamma,.)]\;\vert
\varphi\vert_\infty D(x, y).$$
 By Proposition \ref{contraction} and the mean value relation 
(\ref{TAF1}), there exist $C>0$ and $0\leq r<1$
 such that
 $D( \gamma\cdot x,\gamma\cdot y)\leq C r^k D(x, y)$ whenever $x, 
y\in   \tilde \Lambda_j$, $j\not=l_\gamma$. This 
 leads to the following inequality
\begin{equation}\label{DFr}
[\mathcal L_s^k\varphi] \leq  r_k 
[\varphi] + R_k
|\varphi|_{\infty}
\end{equation}
 where  $r_k = \Bigl(C r^k\Bigr) \;  |\mathcal L_{s}^k1|_{\infty}$ and $R_k =
\sum_{\gamma \in \Gamma(k)} [   w_s(\gamma,.)]$.
Observe that
$$ \limsup_k r_k^{1/k}=
r \limsup_k |\mathcal L_{s}^k1|_{\infty}^{1/k} =
r \rho_s(\infty)$$ where  $\rho_s(\infty)$ is the  spectral 
radius of
the positive operator  $\mathcal L_{s}$ on $C( \tilde \Lambda(\Gamma))$.
Inequality  (\ref{DFr}) is crucial in the Ionescu-Tulcea-Marinescu 
theorem for
quasi-compact operators. By Hennion's work \cite{H}, it implies 
that the essential
spectral radius of $\mathcal L_s$ on ${\bf Lip}( \tilde \Lambda) $ is less than $ 
r \rho_s(\infty)$ ; in
other words, any spectral value of $\mathcal L_s$ with modulus strictly larger 
than $
r \rho_s(\infty)$
is an eigenvalue with finite multiplicity and is  isolated in the spectrum  of
$\mathcal L_s$.

This implies
 in particular  $ \rho_s =  \rho_s(\infty)$. Indeed, the 
inequality
 $ \rho_s\geq   \rho_s(\infty)$ is obvious since the 
function $1$ 
 belongs to ${\bf Lip}( \tilde \Lambda)$.
Conversely, the strict inequality  
would imply   the existence of  a function $\phi \in  {\bf Lip}( \tilde \Lambda)$ such that ${\mathcal 
L}_{s}\phi = \lambda \phi$ for some $\lambda \in \mathbb C$ of modulus 
$>  \rho_s(\infty)$ ; this  yields  $| \lambda| |\phi| \leq
{\mathcal 
L}_{s} |\phi|$ so that $|\lambda |\leq \rho_s(\infty)$. 
Contradiction.

It remains to  control the value $ \rho_s$ in 
the spectrum of ${\mathcal L}_{s}$.  By the above, we know that $ \rho_s$ is an  eigenvalue of $\mathcal L_s$ with  (at least)  one associated 
eigenfunction 
$h_s \geq 0$. This function is strictly positive on $\tilde \Lambda$: otherwise, there  exist $1\leq j \leq p+q$ and  a point  $y_0\in  \tilde \Lambda_j$ such that  $h_s(y_0)=0$.
The equality ${\mathcal L}_{s}h_s(y_{0})= 
 \rho_sh_s(y_{0})$ implies 
  $h_s(\gamma \cdot y_{0}) = 0$ for any $\gamma \in \Gamma$ with 
last letter  $\neq j$.
The minimality of the action of $\Gamma$ 
on $\Lambda$ and the fact that $\Gamma \cdot x_0$ accumulates on $\Lambda$ implies $h_s= 0$ on $\tilde \Lambda$.  Contradiction.

In order to prove that  $ \rho_s$ is a simple eigenvalue of $\mathcal L_s$ on  $ {\bf Lip}( \tilde \Lambda)$, we use a classical argument in probability theory related to the ''Doob transform'' of a sub-markovian transition operator. 
 For any $s \geq \delta$, we denote $P_s$ the operator defined formally   by: for any  bounded Borel function $\phi:  \tilde \Lambda \to \mathbb C$ and $x \in   \tilde \Lambda$,
$$
P_s \phi(x)= {1\over \rho  h_s (x)}\mathcal L (h_s\phi)(x)={1\over \rho  h_s (x)}\sum_{\gamma \in \Gamma(1)}e^{-\delta \tilde b(\gamma, x)} h (\gamma\cdot x)
\phi(\gamma\cdot x).
$$
The  iterates  of   $P_s $ are given by: $P_s ^0= {\rm Id}$ and for $k\geq 1$
\begin{equation}\label{iteresPs}
P_s ^k\phi(x)= \int_X\phi(y)P_s ^k(x, dy)= {1\over \rho_s ^k  h_s  (x)}\sum_{\gamma \in \Gamma(k)}e^{-\delta b  (\gamma, x)}h  (\gamma\cdot x)\phi(\gamma\cdot x).w
\end{equation}

The operator $P_s $ acts on ${\bf Lip}( \tilde \Lambda) $  as  a Markov operator, i.e.  $P_s \phi\geq 0$ if $\phi  \geq 0$ and $P_s {\bf 1} = {\bf 1}$. It inherits the spectral properties of $\mathcal L_s$ and is in particular 
quasi-compact  with   essential spectral radius 
 $<1$. The spectral value $1$ is   an eigenvalue and it remains to prove that the associated eigenspace is $\mathbb C \cdot 1$. Let $f \in {\bf Lip}( \tilde \Lambda)$
such that $P_sf=f$ and $1\leq j \leq p+q$ and $y_0\in  \tilde \Lambda_j$  such that $\vert f(y_{0})\vert = \vert f\vert_{\infty}$. 
An argument of convexity applied to 
the inequality  $P|f|\leq |f|$ readily implies 
$|f(y_{0})|= |f(\gamma\cdot y_{0})|$ for any $\gamma \in \Gamma$ 
with last letter $\neq j$; by minimality of the action of
$\Gamma$ on $ \tilde \Lambda$, it follows that  the modulus of $f$ is
constant  
on $\tilde \Lambda$. Applying  again an argument of convexity, 
the
minimality of the action of
$\Gamma$ on $ \tilde \Lambda$ and the fact that $\Gamma\cdot x_0$ accumulates on  $\Lambda$, one 
proves that $f$ is  in fact constant on $ \tilde \Lambda$. Finally, the eigenspace  of $\mathcal L_s$ associated with $ \rho_s$ equals $\mathbb C \cdot 1$.

Similarly, using the fact that $\ell \geq 3$, we may prove that the  peripherical spectrum of $\mathcal L_s$, i.e. the eigenvalues $\lambda$ with $|\lambda| =  \rho_s$, is reduced to $\rho_s$; we refer the reader to Proposition III.4 of \cite{BP} and Proposition 8.6 of \cite{V}.
 \rightline{$\Box$}

Expression (\ref{iteresPs}) yields to the following.

\begin{notations} For any $s\geq \delta,$ any $x \in   \tilde \Lambda,  $ any $ k\geq 0$ and  any  $\gamma \in \Gamma(k)$,  set
\begin{eqnarray}\label{poids-proba}
p_s (\gamma, x)&:=& {1\over \rho_s ^k }{h_s  (\gamma\cdot x)\over  h_s  (x)} w_s(\gamma, x).
\end{eqnarray}
\end{notations}
As for the $w_s(\gamma, \cdot)$,  these ``weight functions'' are positive  and satisfy   the cocycle property 
$$ 
p_s (\gamma_1\gamma_2, x)= p_s (\gamma_1, \gamma_2\cdot x)\cdot p_s (\gamma_2, x)
$$
for any $s\geq \delta,  x \in   \tilde \Lambda$ and $\gamma_1, \gamma_2 \in \Gamma$.
 Let us emphasize  that $\displaystyle \sum_{\gamma\in \Gamma(k)} p_s (\gamma, x)=1$; in other words, the operator $P_s$ is markovian. 

\rightline{$\Box$}

\begin{coro} \label{coroconvergent}
The group $\Gamma$ is convergent if and only if  $\rho_\delta<1$.
\end{coro}
Proof.    If  $\rho_\delta=\rho_\delta(\infty) <1$ then  $\rho_s<1$ for any $s\geq\delta$, since  $s\mapsto \rho_s(\infty)=\rho_s$ is decreasing on $[\delta, +\infty[$.
Equality  (\ref{criticalexponent}) implies   $\delta_\Gamma \leq  \delta$ and so $\delta_\Gamma=\delta$; by (\ref{divergence-convergence}), it follows  
  that $\Gamma$  is convergent.

Assume now $\rho_\delta\geq1$. When $\Gamma$ is non exotic,  it is divergent by \cite{DOP}. Otherwise, $\delta_\Gamma=\delta$ and  since the eigenfunction $h_\delta$ is non negative on $ \tilde \Lambda$, we have, for any $k\geq 1$ and $x \in  \tilde \Lambda$
$$
\mathcal L_\delta ^k1(x)\asymp
\mathcal L_\delta ^kh_\delta(x)=\rho_\delta^k h_\delta(x)\asymp \rho_\delta^k.
$$
Consequently $\displaystyle \sum_{k\geq 0} \mathcal L_\delta ^k1(x)=+\infty$ and the group $\Gamma$  is divergent, by (\ref{divergence-convergence}).
 
 \rightline{$\Box$}


\section{Counting for convergent groups}\label{sec:Counting}

Throughout this section we assume that $\Gamma$ is convergent on $(X, g )$; by Corollary \ref{coroconvergent} it is equivalent to the fact that $\rho_\delta <1$.

For any $\phi \in {\bf Lip}( \tilde \Lambda)$, any $x \in   \tilde \Lambda$ and $R>0$, let us denote by $M(R, \phi\times\cdot\ )(x)$ the measure on $\mathbb R$ defined by: 
 $$
 M(R, \phi\otimes u)(x):= \sum_{\gamma \in \Gamma} e^{-\delta \tilde b(\gamma, x)}\phi(\gamma\cdot x)u(-R+\tilde b(\gamma, x)).
 $$
It holds  $0\leq M(R, \phi\otimes u)(x) <+\infty$ when $u$ has a compact support in $ \mathbb R $ since the group $\Gamma$ is discrete.
  
 The orbital function of $\Gamma$ may be decomposed as
 $$N_\Gamma(R)= e^{\delta R} \sum_{n\geq 0}M(R, {\bf 1} \otimes e_n)(x_0)
 $$ with $e_n(t):= e^{\delta  t}{\bf 1}_{]-(n+1), -n]}(t).
 $
Hence, Theorem \ref{theo:Comptage} is a direct consequence of the following statement.

 \begin{prop} \label{convergence vague}
For any positive function $\phi \in {\bf Lip}( \tilde \Lambda)$ and any $x \in   \tilde \Lambda$, there exists $C_\phi(x)>0$ such that for any continuous function  $u: \mathbb R\to \mathbb R$ with compact support,
$$
 \lim_{R\to +\infty}{R^{\alpha}\over L(R)}M(R, \phi\otimes u)(x) = C_\phi(x) \int_{\mathbb R}u(t) {\rm d}t.
$$
\end{prop}
%
%
This section is devoted to the proof of Proposition \ref{convergence vague}. From now on, we fix a positive function $\phi\in {\bf Lip}( \tilde \Lambda)  $ and  a continuous function $u: \mathbb R\to \mathbb R^+$ with compact support.
Let us decompose    $M(R, \phi\otimes u)(x)$ as 
$$
M(R, \phi\otimes u)(x)=\sum_{k \geq 0} M_k(R, \phi\otimes u)(x)
$$
with
$$
 M_k(R, \phi\otimes u)(x):= 
 \sum_{\gamma \in \Gamma(k)} e^{-\delta \tilde b(\gamma, x)}\phi(\gamma\cdot x) u(-R+\tilde b(\gamma, x)).
$$
Thus, it is natural to  associate to $P_s, s\geq \delta, $ a new  transition operator $\widetilde P_s $ on $  \tilde \Lambda \times \mathbb R$, setting: for any $\phi\in {\bf Lip}( \tilde \Lambda)  $,  any Borel function $v: \mathbb R\to \mathbb R$ and any $(x, s) \in   \tilde \Lambda)\times \mathbb R,$
\begin{eqnarray*}\label{opeP}
\widetilde P_s (\phi\otimes v)(x, t)&=& {1\over \rho   h_s (x)}\sum_{\gamma \in \Gamma(1)}e^{-s \tilde b(\gamma, x)} h_s (\gamma\cdot x)\phi(\gamma\cdot x) u(t+\tilde b(\gamma, x))\\
&=& \sum_{\gamma \in \Gamma(1)}p_s(\gamma, x) \phi(\gamma\cdot x) u(t+\tilde b(\gamma, x))\notag
\end{eqnarray*}
Notice that $\widetilde P_s $ is a also a Markov operator on $  \tilde \Lambda \times \mathbb R$; it commutes with the action of translations on $\mathbb R$ and  one usually says that it defines a semi-markovian  random walk on $  \tilde \Lambda \times \mathbb R$. 
Its iterates  are given by: 
$\widetilde P_s^0= {\rm Id}$ and, for any $k\geq 1$,
$$
\widetilde P_s ^k(\phi\otimes v)(x, s)= \sum_{\gamma \in \Gamma(k)} p_s\gamma, x)\phi(\gamma\cdot x)u(s+\tilde b(\gamma, x)).
$$
From now on, to lighten notations we  write $P = P_\delta$, $\tilde P = \tilde P_\delta$, $h = h_\delta$, $  p =  p_\delta$ and $\rho = \rho_\delta<1$.
We rewrite the quantity $ M_k(R, \phi\otimes u)(x)$ as
$$
 M_k(R, \phi\otimes u)(x)=\rho ^kh (x) \widetilde P^k\left( {\phi\over h }\otimes u\right)(x, -R),
$$
so that, 
\begin{equation}\label{Maspotential}
M(R, \phi\otimes u)(x)=h (x)\sum_{k\geq 0}\rho ^k \widetilde P^k\left({\phi\over h }\otimes u\right)(x, -R).
\end{equation}
We first control the   behavior as $R\to +\infty$  of the quantity $ M_1(R, \phi\otimes u)(x)$.

 \begin{prop}\label{prop:AsymPk1}
 For any continuous function $u: \mathbb R\to \mathbb R$ with compact support, there exists a constant $C_u>0$ such that, for any $\varphi \in  {\bf Lip}( \tilde \Lambda)$,   any $ x  \in \tilde  \Lambda$  and $R\geq 1$,

\begin{equation}\label{majorforbusemann}
 \Big\vert \widetilde P (\varphi\otimes u)(x, -R) \Big\vert \leq C_u  \Vert \varphi \Vert_\infty   \times  {  L(R)\over R^{\alpha}}.
\end{equation}
Furthermore,
\begin{equation}\label{asymptforbusemann}
 \lim_{R\to +\infty} {R^{\alpha} \over L(R)}
 \widetilde P (\varphi\otimes u)(x, -R)=\sum_{j=1}^p C_j(x) \varphi(x_j)\int_{\mathbb R}u(t) {\rm d}t, 
 \end{equation}
where $C_j$ is defined by: for  $1\leq j \leq p$,
\begin{equation}\label{constanteC}
 C_j(x):= c_j{h (x_j)  \over \rho  h (x)} \times \left\{ \begin{array}{cll}
\displaystyle  e^{2 \delta(x_j\mid x)_{\bf o}}& {\rm when} &x\in  \Lambda \backslash \tilde \Lambda_j;
 \\
e^{\mathcal B_{x_j}(o, g\cdot o) + d(o, g\cdot o)}& {\rm when} &x = g\cdot x_0 \notin \tilde \Lambda_j;\\
 
 \ & \ 
 \\
 0
& {\rm otherwise.}\ \ &\ 
\end{array}\right.
\end{equation}
\end{prop}

\noindent Proof. 
 Let $x\in \tilde \Lambda$ be fixed and assume that the support of $u$ is included in the interval $[a, b]$. For any $R\geq -a$, it holds
$$
 \widetilde P (\varphi\otimes u)(x, -R) = 
\frac{1}{\rho h(x)}\sum_{j = 1}^{p+q}\sum_{\gamma \in \Gamma_j}e^{-\delta \tilde b(\gamma  x)}{\bf 1}_{x\notin \tilde \Lambda_j}h(\gamma\cdot x) \varphi(\gamma\cdot x) u(-R+\tilde b(\gamma, x) ).
$$
It follows from hypotheses $  {\bf H_2}$ and $  {\bf H_3}$     and Fact \ref{lienentrebetd} that for any $j = 1,..., p+q$, there exists a constant   $ K_j>0$  such that for any $R\geq 1$,
$$\sum_{\stackrel{\gamma \in \Gamma_j}{R+a \leq \tilde b(\gamma, x) \leq R+b}} e^{-\delta \tilde b(\gamma, x)} \leq  K_j (b-a) \frac{L(R)}{R^\alpha}.$$
Together with the fact that $L$ has slow variation, this implies (\ref{majorforbusemann}).

Now, in order to establish (\ref{asymptforbusemann}), it is sufficient to prove  that for  any $j = 1,..., p+q$,
\begin{equation}\label{asymptforbusemanntermeparterme}
 \lim_{R\to +\infty} {R^{\alpha} \over L(R)}
 \sum_{\gamma\in \Gamma_j}   p(\gamma, x)\varphi(\gamma\cdot x) u(-R + \tilde b(\gamma, x)) = C_j(x) \varphi(x_j)\int_{\mathbb R}u(t) {\rm d}t,
\end{equation}
where $C_j(x)$ is given by (\ref{constanteC}) for $1\leq j \leq p$  and  $C_j(x) = 0$  for $j = p+1,..., q$.   By a  classical approximation argument, we may assume that 
 $u$ is the characteristic function of the  interval $[a, b]$; it yields
$$ 
\sum_{\gamma\in \Gamma_j}  p(\gamma, x)\varphi(\gamma\cdot x) u(-R + \tilde b(\gamma, x))  = 
\frac{1}{h(x)} \sum_{\stackrel{\gamma \in \Gamma_j}{R+a \leq \tilde b(\gamma,   x) \leq R+b}}e^{-\delta \tilde b(\gamma  x)}{\bf 1}_{x\notin \tilde \Lambda_j}h(\gamma\cdot x) \varphi(\gamma\cdot x).
$$
First, assume that  $x = g\cdot x_0$  belongs to  $\Gamma\cdot x_0$.   For  any $j = 1,..., p$ and  any  $\gamma\neq Id$ in  $ \Gamma_j$, the sequence 
$( \gamma^n \cdot o)_{n \geq 1}
$ tends to  $x_j$ as $n\to \pm\infty$;   it yields
\begin{eqnarray*}
 \tilde b(\gamma^n, x) - d(o, \gamma^n \cdot o) 
 &=& d(\gamma^{-n}\cdot o, g\cdot o) - d(\gamma^{-n}\cdot o, o) - d(o, g\cdot o)\\
 &{\stackrel{n \to \pm \infty}{\longrightarrow}}& - \mathcal B_{x_j}(o, g\cdot o) - d(o, g\cdot o).
\end{eqnarray*}
 When $x\in \Lambda$,    Fact \ref{lienentrebetd} yields 
$$
 \lim_{n \to \pm \infty} \tilde b(\gamma^n, x) - d(o, \gamma^n \cdot o)= -2(x_j\mid x).
 $$
Eventually, by  hypotheses   $  {\bf H_2}$  and $  {\bf H_3}$, for any $1\leq j \leq p+q$,    
$$
 \lim_{R\to +\infty} {R^{\alpha} \over L(R)}
 \sum_{\stackrel{\gamma\in \Gamma_j}{R+a \leq d(o, \gamma \cdot o) \leq R+b}} \tilde p(\gamma, x) = C_j(x) |b-a|.
$$
Hence,
$$
 \lim_{R\to +\infty} {R^{\alpha} \over L(R)}
 \sum_{\gamma\in \Gamma_j} \tilde p(\gamma, x)\varphi(\gamma\cdot x) u(-R + \tilde b(\gamma, x)) 
= C_j(x) \varphi(x_j) \vert b-a\vert.
$$

\rightline{$\Box$}

%

Now, we extend (\ref{majorforbusemann})  and (\ref{asymptforbusemann}) to the powers $\widetilde P^k, k \geq 1,$ of the Markov operator  $\widetilde P$.

\begin{prop}\label{mkmaj} 
For any  continuous function  $u:\mathbb R \to \mathbb R^+$ with compact support, there exists a constant $C_u>0$ such that,  for any $\varphi \in {\bf Lip}(\tilde \Lambda)$,   any $x \in \tilde \Lambda$, any $k \geq 1$  and any $R\geq 1$,
\begin{equation} \label{mkmaj-formule}
\Bigl\vert \widetilde P^k\left(\varphi\otimes u\right)(x, -R)\Bigr\vert \leq  C_u  \  k^2 \ \Vert \varphi \Vert_\infty   \times {L(R)\over R^{\alpha}}.
\end{equation}
\end{prop}

\begin{prop}\label{mkequi}
For any   continuous function  $u:\mathbb R \to \mathbb R^+$ with compact support,  any $\varphi \in {\bf Lip}(\tilde \Lambda)$,  any $x \in   \tilde \Lambda$ and any $k \geq 1,$
\begin{equation} \label{mkequi-formule}
\lim_{R\to +\infty}
{ R^{\alpha}\over L(R)}  \widetilde P^k\left(\varphi\otimes u\right)(x, -R)= \sum_{j = 1}^p \left(\sum_{l=0}^{k-1}P^l C_j(x)P^{k-1-l}\varphi(x_j)\right) \int_{\mathbb R} u(t) {\rm d}t
\end{equation}
where, for any $1\leq j\leq p$, the Lipschitz functions is $C_j: \tilde \Lambda \to \mathbb R$ is   given by (\ref{constanteC}).
\end{prop}

Proposition \ref{convergence vague} follows immediately from these statements and (\ref{Maspotential}). Indeed, Propositions \ref{mkmaj} and       \ref{mkequi} and the dominated convergence theorem yield   
$$
\lim_{R\to +\infty} \frac{R^\alpha}{L(R)} M(R, \phi\otimes u)(x) = \left( h (x)\sum_{k\geq 1}\rho ^k \sum_{j = 1}^p\left(\sum_{l=0}^{k-1}P^l C_j(x) P^{k-1-l}\left(\frac{\phi}{h}\right)(x_j)\right)\right)\times \int_{\mathbb R} u(t) {\rm d}t.
$$

\rightline{$\Box$}

Let us now prove Propositions \ref{mkmaj} and \ref{mkequi}. 
For the convenience of the reader, we  assume that  all subgroups  $\Gamma_j, 1\leq j \leq p+q,$  are parabolic. Hence, they  have a unique fixed point at infinity $x_j$ and for any $x\in \tilde \Lambda$, it holds
$$
\lim_{\stackrel{\gamma\in \Gamma_j}{d(o, \gamma\cdot o) \to +\infty}} \gamma\cdot x = x_j.
$$ 
Namely, if  one of the non-influent elementary group $\Gamma_j, p+1\leq j \leq p+q,$ was  generated by some hyperbolic isometry $h_j$, we would have in the next proofs to distinguish between positive and negative power of $h_j$ and  this would only overcharge our notations without interest.

\noindent Proof of Proposition \ref{mkmaj}. We apply here overestimations given in \cite{V}, whose proofs follow  the approach developed in  \cite{G}. We set $\alpha = 1+\beta$ with $0< \beta <1$; this restriction on the values of the parameter $\beta$ is of major importance to get the following estimations. Following \cite{V}, we introduce the non negative sequence $(a_k)_{k \geq 1}$ defined implicitely by $\displaystyle {a_k^\beta\over L(a_k)}=k$ for any $k \geq 1$. By Propositions A.1 and A.2 in \cite{V}, there exists a constant $C_1=C_1(u)>0$ such that,  for any $\varphi \in {\bf Lip}(\tilde \Lambda)$,   any $x \in \tilde \Lambda$, any $k \geq 1$ and any $R\geq  1$,

$\bullet$  if  $1\leq R\leq 2a_k$ then 
$\quad \displaystyle 
\Bigl\vert \widetilde P^k\left(\varphi\otimes u\right)(x, -R)\Bigr\vert \leq  C_1   \Vert \varphi \Vert_\infty   \times{1\over a_k};
$

$\bullet$ if  $R\geq 2a_k$ then 
 $\displaystyle 
\Bigl\vert \widetilde P^k\left(\varphi\otimes u\right)(x, -R)\Bigr\vert \leq   C_1 k   \Vert \varphi \Vert_\infty   \times {L(R)\over R^{1+\beta}}. 
$

 \noindent  The definition of the $a_k$ yields, for $1\leq R\leq 2a_k$, 
$$
{1\over a_k}=k{L(a_k)\over a_k^{1+\beta}}\leq {k\over 2^{1+\beta}}\times {L(R)\over R^{1+\beta}}\times {L(a_k)\over L(R)}.$$
By Potter's lemma (see \cite{V}, lemma 3.4),  it  exists $C_2>0$ such that $\displaystyle {1\over a_k}\leq C_2 k^2\times {L(R)\over R^{1+\beta}}
$ for   $R\geq 1$ great enough. We set $C= \max( C_1, C_2).$

\rightline{$\Box$}

\noindent Proof of Proposition \ref{mkequi}.
We  work by induction. By Proposition \ref{prop:AsymPk1}, convergence (\ref{mkequi-formule}) holds for $k=1$. Now, we assume that it holds for some $k\geq 1$.
Let $R>0$ and $r\in [0, R/2]$ be fixed. 
Recall that 
\begin{eqnarray*}
\widetilde P^{k+1}\left(\varphi\otimes u\right)(x, -R)&=&  
\sum_{\gamma \in \Gamma(k+1)}p (\gamma, x) \varphi(\gamma\cdot x) u(-R+\tilde b(\gamma, x))\\
&=&
\sum_{\gamma \in \Gamma(k)} \sum_{\beta \in \Gamma(1)}p (\gamma, \beta\cdot x) 
p (\beta,  x)\varphi(\gamma \beta\cdot x) u\Bigl(-R+\tilde b(\gamma, \beta\cdot x)+
\tilde b(\beta,   x)\Bigr).
\end{eqnarray*}
We decompose $\widetilde P^{k+1}\left(\varphi\otimes u\right)(x, -R)$ as
$ 
 A_k(x, r, R)+B_k(x, r, R)+C_k(x, r, R)
 $
where
\begin{eqnarray*}
A_k(x, r, R) &:=&  \sum_{\gamma \in \Gamma(k)} \sum_{\stackrel{\beta \in \Gamma(1)}{d({\bf o}, \beta\cdot{\bf o})\leq r}}
p (\gamma, \beta\cdot x) 
p (\beta,x)\varphi(\gamma \beta\cdot x) u\Bigl(-R+\tilde b(\gamma, \beta\cdot x)+
\tilde b(\beta,   x)\Bigr),\\
B_k(x, r, R) &:=& \sum_{\stackrel{\gamma \in \Gamma(k)}{d({\bf o}, \gamma\cdot{\bf o})\leq r}}
\sum_{\stackrel{\beta \in \Gamma(1)}{d({\bf o}, \beta\cdot{\bf o})>r}}
p (\gamma, \beta\cdot x) 
p (\beta,  x)\varphi(\gamma \beta\cdot x) u\Bigl(-R+\tilde b(\gamma, \beta\cdot x)+
\tilde b(\beta,   x)\Bigr)\\
{\rm and} \ \ 
C_k(x, r, R) &:=&
 \sum_{\stackrel{\gamma \in \Gamma(k)}{ d({\bf o}, \gamma\cdot{\bf o})>r}}
\sum_{\stackrel{\beta \in \Gamma(1)}{d({\bf o}, \beta\cdot{\bf o})>r}}\ 
 p (\gamma, \beta\cdot x) 
p (\beta,  x)\varphi(\gamma \beta\cdot x) u\Bigl(-R+\tilde b(\gamma, \beta\cdot x)+
\tilde b(\beta,  x)\Bigr).
\end{eqnarray*}

\noindent \underline{ \bf Step 1.} Let us first prove that 
\begin{equation}\label{Akrfixe}
\lim_{R\to +\infty}
{R^{\alpha} \over L(R) } A_k(x,r, R) = 
\sum_{\stackrel{\beta \in \Gamma(1)}{d({\bf o}, \beta\cdot {\bf o})\leq r}}
p(\beta, x)\times \lim_{R\to +\infty}
\ {R^{\alpha} \over L(R) }  \tilde P^k\left(\varphi\otimes u\right)(\beta \cdot x, -R).
\end{equation}
Indeed, 
the set of $\beta \in \Gamma(1)$ such that $d({\bf o}, \beta\cdot{\bf o})\leq r$ is   finite   and  $\tilde b(\beta,  x)\leq r$ for such an isometry $\beta$; furthermore,   if $p(\beta, x)\neq 0$ then 
 ${R\over 2}  \leq R-\tilde b(\gamma  \beta\cdot x)\leq R+C$ where $C>0$ is the constant which appears in  Property \ref{triangle}. Using the induction hypothesis,  it yields, for   any   $\beta \in \Gamma(1)$  such that $d({\bf o}, \beta\cdot{\bf o})\leq r$,
  
\noindent $\displaystyle 
\lim_{R\to +\infty} {R^{\alpha}\over L(R)} p (\beta,  x)
\sum_{\gamma \in \Gamma(k)}
p (\gamma, \beta\cdot x) 
\varphi(\gamma \beta\cdot x) u\Bigl(-R+\tilde b(\beta, x)+
\tilde b(\gamma, \beta \cdot x)\Bigr)$

$ \displaystyle \qquad \qquad \qquad \qquad\qquad \qquad \qquad \qquad =
p(\beta, x)\times \lim_{R\to +\infty}
\ {R^{\alpha} \over L(R) } \tilde P^k\left(\varphi\otimes u\right)(\beta \cdot x,R).
$

\noindent Convergence (\ref{Akrfixe}) follows, 
summing over $\beta$.
It yields
\begin{equation}\label{Akrinfty}
\lim_{r\to +\infty}
\lim_{R\to +\infty}
\ {R^{\alpha} \over L(R) } A_k(x, r, R)= \sum_{j = 1}^p \left(\sum_{l=1}^{k}P^lC_j(x)P^{k-l}\varphi(x_j)\right)
\times \int_{\mathbb R}u(t) {\rm d}t.
\end{equation}

 \noindent \underline{\bf Step 2.}
We prove that there 
  exists $\epsilon(r)>0$, with $\displaystyle \lim_{r\to +\infty}\epsilon(r) = 0$, such that, for any $k\geq 1$,
\begin{eqnarray}\label{Bkrfixe}
\liminf_{R\to +\infty}
\ {R^{\alpha} \over L(R) }B_k(x,r, R) &{\stackrel{\epsilon(r)}{ \simeq }}& 
\limsup_{R\to +\infty}
\ {R^{\alpha} \over L(R) }B_k(x,r, R) \notag  \\
&{\stackrel{\epsilon(r)}{ \simeq }}&
\sum_{j = 1}^p \sum_{\stackrel{\gamma\in \Gamma(k)}{d(o, \gamma\cdot o)\leq r}}  p(\gamma, x_j) \varphi(\gamma\cdot x_j)C_j(x) \int_{\mathbb R} u(t) {\rm d}t,
\end{eqnarray}
where we write $a\  {\stackrel{\epsilon }{ \simeq }} \ b$  if $\displaystyle 1-\epsilon \leq \frac{a}{b} \leq 1+\epsilon $.
Since each $\Gamma_j$ has a unique fixed point, there exists a map $\epsilon : (0, +\infty) \to (0, +\infty)$ which tends to $0$ as $r \to +\infty$,   such that 
$$
 \frac{p(\gamma, \beta\cdot x)}{p(\gamma,  x_j)}  \ {\stackrel{\epsilon(r) }{ \simeq }} \   \ 1 
$$
for any $j = 1,..., p+q$, any $\beta\in \Gamma_j$ with $d(o, \beta\cdot o)\geq r$, any $x\in \tilde \Lambda$ and any $\gamma\in \Gamma$ with $l_\gamma \neq j$.

The set of $\gamma \in \Gamma(k)$ such that $d({\bf o}, \gamma\cdot{\bf o})\leq r$ is   a finite subset of $\Gamma(k)$; furthermore, for such  $\gamma$ and any $\beta \in \Gamma(1)$, it holds    ${R\over 2}\leq R  -\tilde b(\gamma, \beta\cdot x) \leq R+C, $ as above.
Therefore,

\vspace{3mm}

\noindent $
\displaystyle \sum_{\stackrel{ \gamma \in \Gamma(k)}{d({\bf o}, \gamma\cdot{\bf o})\leq r}}
\sum_{
\stackrel{ \beta \in \Gamma(1)}{d({\bf o}, \beta\cdot{\bf o})>r}} 
 p (\gamma, \beta\cdot x) 
p (\beta,  x)\varphi(\gamma \beta\cdot x) u\Bigl(-R+\tilde b(\gamma, \beta\cdot x)+
\tilde b(\beta,  x)\Bigr) $
\\
\indent 
$\qquad \qquad  
\displaystyle  \ {\stackrel{\epsilon(r) }{ \simeq }} \  
\sum_{j = 1}^{p+q} 
\sum_{\stackrel{ \gamma \in \Gamma(k)}{d({\bf o}, \gamma\cdot{\bf o})\leq r}}
 p (\gamma,  x_j) \varphi(\gamma \cdot  x_j)
\sum_{
\stackrel{ \beta \in \Gamma_j}{d({\bf o}, \beta\cdot{\bf o})>r}} 
p (\beta,  x) u\Bigl(-R+\tilde b(\gamma, \beta\cdot x)+
\tilde b(\beta,  x)\Bigr) 
$

\vspace{3mm}

\noindent Convergence (\ref{Bkrfixe}) follows, using  (\ref{asymptforbusemanntermeparterme}).
In particular, letting $r \to +\infty$, it holds
\begin{eqnarray}  \label{Bkrinfty}
\lim_{r\to +\infty}
\liminf_{R\to +\infty}
\ {R^{\alpha} \over L(R) }B_k(x,r, R) &=& \lim_{r\to +\infty}
\limsup_{R\to +\infty}
\ {R^{\alpha} \over L(R) }B_k(x,r, R) \notag  \\
&=&
\sum_{j = 1}^{p} P^k\varphi( x_j)C_j(x) \int_{\mathbb R} u(t) {\rm d}t.
\end{eqnarray}

\noindent \underline{\bf Step 3.} We prove  that 
there exists  a constant $C>0$  such that,  for any $R\geq 2r\geq 1$,
\begin{equation}  \label{Ckrfixe}
 C_k(x,r,R)\leq C  k^2  
\Vert \varphi \Vert_\infty 
{L(R)\over R^{\alpha}}
\sum_{n=[r]}^{+\infty} 
{L(n)\over n^{\alpha}}.
\end{equation}

By property \ref{triangle}, the condition $u\Bigl(-R+\tilde b(\gamma  \beta\cdot x)+
\tilde b(\beta,  x)\Bigr)\neq 0$ implies
$$d({\bf o}, \gamma\cdot{\bf o})+d({\bf o}, \beta\cdot{\bf o}) =R\pm c\quad {\rm and} \quad \tilde b(\gamma  \beta\cdot x)+
\tilde b(\beta,  x) =R\pm c \quad ^(\footnote{the notation $A=B\pm c$ means $\vert A-B\vert \leq c$.}^) $$ 
for some constant $c>0$ which depends on $u$.

We decompose  $C_k(x, r, R)$ into $C_k(x, r, R)=C_{k,1}(x, r, R)+ C_{k,2}(x, r, R)$ with 
$$
C_{k,1}(x, r, R):= \sum_{\stackrel{\gamma \in \Gamma(k)}{r<d({\bf o}, \gamma\cdot {\bf o})\leq  R/2}}
\sum_{\stackrel{\beta \in \Gamma(1)}{d({\bf o}, \beta\cdot{\bf o})>r}}\ 
p (\gamma, \beta\cdot x) 
p (\beta,  x)\varphi(\gamma \beta\cdot x) u\Bigl(-R+\tilde b(\gamma, \beta\cdot x)+
\tilde b(\beta,  x)\Bigr).
$$
and
$$
C_{k,2}(x, r, R):= \sum_{\stackrel{\gamma \in \Gamma(k)}{ d({\bf o}, \gamma\cdot {\bf o})\geq R/2}}
\sum_{\stackrel{\beta \in \Gamma(k)}{d({\bf o}, \beta\cdot{\bf o})>r}}\ 
p (\gamma, \beta\cdot x) 
p (\beta,  x)\varphi(\gamma \beta\cdot x) u\Bigl(-R+\tilde b(\gamma, \beta\cdot x)+
\tilde b(\beta,  x)\Bigr).
$$
We control the term $C_{k,1}(x, r, R)$. Assuming $c \geq 1$, one may write
\begin{eqnarray*}
C_{k,1}(x, r, R)&\leq& \Vert \varphi \Vert_\infty 
\Vert u \Vert_\infty
\sum_{n=[r]}^{[R/2]} 
\ 
\sum_{\stackrel{\gamma \in \Gamma(k)}{ d({\bf o}, \gamma\cdot {\bf o})= n\pm c}}
\ 
\sum_{\stackrel{\beta \in \Gamma(k)}{d({\bf o}, \beta\cdot{\bf o})= R-n\pm c}}\ 
p (\gamma, \beta\cdot x) 
p (\beta,  x)\\ 
&\leq &
\Vert \varphi \Vert_\infty 
\Vert u \Vert_\infty
\sum_{n=[r]}^{[R/2]} 
\sum_{\stackrel{\beta \in \Gamma(1)}{d({\bf o}, \beta\cdot{\bf o})= R-n\pm c}}
p (\beta,  x)
\left(
\sum_{\stackrel{\gamma \in \Gamma(k)}{ d({\bf o}, \gamma\cdot {\bf o})=n\pm c}} 
p (\gamma, \beta\cdot x) \right).
\end{eqnarray*}
Using (\ref{mkmaj-formule}),  this  yields, for some constant $C>0$, 
\begin{eqnarray*}
C_{k,1}(x,r,R)&\leq &
C  k^2
\Vert \varphi \Vert_\infty 
\Vert u \Vert_\infty
\sum_{n=[r]}^{[R/2]} 
{L(R-n)\over (R-n)^{\alpha}}
{L(n)\over n^{\alpha}}  
\\
&\leq &
C\ {k^2}
\Vert \varphi \Vert_\infty 
\Vert u \Vert_\infty
{L(R)\over R^{\alpha}}
\sum_{n=[r]}^{+\infty} 
{L(n)\over n^{\alpha}}, 
\end{eqnarray*}
where the last inequality  is based  on the  facts that  $R-n\geq R/2-1$ and $L$ is slowly varying.
The same inequality holds for $C_{k,2}(x,r,R) $, by reversing  in the previous argument the role of $\gamma$ and $\beta$.  Hence, 
\begin{equation}  \label{Ckrinfty}
 \lim_{r\to +\infty}\limsup_{R\to +\infty}\ {R^{\alpha} \over L(R) } C_k(x, r, R)  = 0.
\end{equation}
%
 Proposition \ref{mkequi} follows, combining  (\ref{Akrinfty}),  (\ref{Bkrinfty}) and  (\ref{Ckrinfty}).
 
\rightline{$\Box$}

\end{document}